\documentclass[a4paper]{amsart}    % option draft shows overfull boxes
\usepackage{multirow}
\usepackage{amsmath}
\usepackage{amsfonts}
\usepackage{amssymb}
\usepackage{latexsym}

\renewcommand{\d}{\,{\rm d}} % d for integrals and differentials

\newcommand{\roundb}[1]{\left(#1\right)} % round brackets
\newcommand{\curlyb}[1]{\left\{#1\right\}} % curly brackets
\newcommand{\setb}[2]{\curlyb{#1 \colon #2}} % set
\newcommand{\triple}[3]{J\left(#1,#2,#3\right)}

\newcommand{\Hproduct}[1]{%
  \roundb{\begin{smallmatrix}#1\end{smallmatrix}}
} %
\newcommand{\LHproduct}[1]{%
  \roundb{\begin{matrix}#1\end{matrix}}
} %

\newcommand{\Nbold}{\boldsymbol{N}}
\newcommand{\Lbold}{\boldsymbol{L}}
\newcommand{\Nmin}{N_{\mathrm{min}}}

\newcommand{\Nmax}{N_{\mathrm{max}}}

\newcommand{\Lmax}{L_{\mathrm{max}}}

\newcommand{\Abs}[1]{\left\lvert #1 \right\rvert}
\newcommand{\abs}[1]{\lvert #1 \rvert}
\newcommand{\bigabs}[1]{\bigl\lvert #1 \bigr\rvert}

\newcommand{\norm}[1]{\left\lVert #1 \right\rVert}
\newcommand{\fixednorm}[1]{\lVert #1 \rVert}
\newcommand{\bignorm}[1]{\bigl\lVert #1 \bigr\rVert}

\newcommand{\R}{\mathbb{R}}

\newcommand{\angles}[1]{\langle #1 \rangle}

\newtheorem{theorem}{Theorem}[section]

\theoremstyle{definition}
\newtheorem{definition}[theorem]{Definition}

\theoremstyle{remark}
\newtheorem{remark}[theorem]{Remark}

\numberwithin{equation}{section}
\title{Atlas of products for wave-Sobolev spaces on $\R^{1+3}$}

\author[P. D'Ancona]{Piero D'Ancona}
\address{Department of Mathematics\\
University of Rome ``La Sapienza''\\
Piazzale Aldo Moro 2\\
I-00185 Rome\\ Italy}
\email{dancona@mat.uniroma1.it}
  
\author[D. Foschi]{Damiano Foschi}
\address{Department of Mathematics\\
University of Ferrara\\
Via Macchiavelli 35\\
I-44100 Ferrara\\ Italy}
\email{damiano.foschi@unife.it}

\author[S. Selberg]{Sigmund Selberg}
\address{Department of Mathematical Sciences\\ Norwegian University of Science and Technology\\ N-7491 Trondheim\\ Norway}
\email{sselberg@math.ntnu.no}
\urladdr{www.math.ntnu.no/~sselberg}

\thanks{This paper was written as part of  the international research program on Nonlinear Partial Differential Equations at the Centre for Advanced Study at the Norwegian Academy of Science and Letters in Oslo during the academic year 2008--09.}

\begin{document}

\begin{abstract}
The wave-Sobolev spaces $H^{s,b}$ are $L^2$-based Sobolev spaces on the Minkowski space-time $\R^{1+n}$, with Fourier weights are adapted to the symbol of the d'Alembertian. They are a standard tool in the study of regularity properties of nonlinear wave equations, and in such applications the need arises for product estimates in these spaces. Unfortunately, it seems that with every new application some estimates come up which have not yet appeared in the literature, and then one has to resort to a set of well-established procedures for proving the missing estimates. To relieve the tedium of having to constantly fill in such gaps ``by hand'', we make here a systematic effort to determine the complete set of estimates in the bilinear case. We determine a set of necessary conditions for a product estimate $H^{s_1,b_1} \cdot H^{s_2,b_2} \hookrightarrow H^{-s_0,-b_0}$ to hold. These conditions define a polyhedron $\Omega$ in the space $\R^6$ of exponents $(s_0,s_1,s_2,b_0,b_1,b_2)$. We then show, in space dimension $n=3$, that all points in the interior of $\Omega$, and all points on the faces minus the edges, give product estimates. We can also allow some but not all points on the edges, but here we do not claim to have the sharp result. The corresponding result for $n=2$ and $n=1$ will be published elsewhere.
\end{abstract}

\maketitle
\tableofcontents

% main text

\section{Wave-Sobolev spaces}\label{A}

Define the Fourier transform of a Schwartz function $u \in \mathcal S(\R^{1+n})$ by
$$
  \widetilde u(\tau,\xi)
  =
  \iint e^{-i(t\tau+x\cdot\xi)} u(t,x) \d t \d x,
$$
where $(t,x)$ and $(\tau,\xi)$ belong to $\R \times \R^n = \R^{1+n}$; $\tau$ and $\xi$ will be called the temporal and spatial frequencies, respectively.

\begin{definition}\label{A:Def1} Given $s,b \in \R$, the \emph{wave-Sobolev space} $H^{s,b}=H^{s,b}(\R^{1+n})$ is the completion of $\mathcal S(\R^{1+n})$ with respect to to the norm
$$
  \norm{u}_{H^{s,b}}
  =
  \norm{ \angles{\xi}^s \angles{\abs{\tau}-\abs{\xi}}^b \widetilde u(\tau,\xi) }_{L^2_{\tau,\xi}},
$$
where $\angles{\cdot} = (1+\abs{\cdot}^2)^{\frac12}$. We shall refer to the weights $\angles{\xi}^s$ and $\angles{\abs{\tau}-\abs{\xi}}^b$ as \emph{elliptic} and \emph{hyperbolic}, respectively.
\end{definition}

By way of comparison, the elliptic weight is a familiar aspect of the standard Sobolev space $H^s=H^s(\R^n)$, obtained as the completion of $\mathcal S(\R^n)$ with respect to the norm $\norm{f}_{H^{s}} = \fixednorm{ \angles{\xi}^s \widehat f(\xi) }_{L^2_{\xi}}$, where $\widehat f(\xi) = \int e^{-i x\cdot\xi} f(x) \d x$.

The hyperbolic weight, on the other hand, reflects the fact that the $H^{s,b}$-norm is adapted to the wave operator, or d'Alembertian, $\square = -\partial_t^2 + \Delta_x$, whose symbol is $\tau^2 - \abs{\xi}^2$.

For details about the history of the wave-Sobolev spaces and applications to nonlinear wave equations, we refer to the survey article \cite{Selberg:2002b}. In the applications, the need frequently arises for product estimates of the form
\begin{equation}\label{A:10}
  H^{s_1,b_1} \cdot H^{s_2,b_2} \hookrightarrow H^{-s_0,-b_0}.
\end{equation}
Explicitly, this means that there exists $C = C(s_0,s_1,s_2,b_0,b_1,b_2; n)$ such that
\begin{equation}\label{A:20}
  \norm{uv}_{H^{-s_0,-b_0}} \le C \norm{u}_{H^{s_1,b_1}}\norm{v}_{H^{s_2,b_2}}
\end{equation}
for all $u,v \in \mathcal S(\R^{1+n})$. 

\begin{definition}\label{A:Def2} If~\eqref{A:20} holds, we say that the exponent matrix
$$
  \LHproduct{s_0 & s_1 & s_2 \\ b_0 & b_1 & b_2}
$$
is a \emph{product}.
\end{definition}

Many product estimates have appeared in the literature (see \cite{Selberg:2002b} for examples and references), but so far no systematic effort has been made to determine necessary and sufficient conditions on $\Hproduct{s_0 & s_1 & s_2 \\ b_0 & b_1 & b_2}$ for it to be a product. The present paper is the result of our efforts to fill this gap. We remark that the utility of these product estimates is not limited to simple products: Many bilinear null form estimates can also be reduced to this form, since the null symbol can be estimated in terms of the weights appearing in the $H^{s,b}$-norm. See, e.g., \cite{Selberg:2007d}.

It turns out that there are 21 necessary conditions of the form
\begin{equation}\label{A:30}
  \sigma_0 s_0 + \sigma_1 s_1 + \sigma_2 s_2 + \beta_0 b_0 + \beta_1 b_1 + \beta_2 b_2 \ge 0.
\end{equation}
Each such condition determines a half-space in the space $\R^6$ of coefficients $\Hproduct{s_0 & s_1 & s_2 \\ b_0 & b_1 & b_2}$. Taken together, these 21 conditions, which are listed in the next section, determine a convex polyhedron $\Omega$ in $\R^6$. The boundary of $\Omega$ consists of \emph{faces}, which are polyhedrons contained in the hyperplanes corresponding to equality in one of the conditions of the form~\eqref{A:30}. The intersection of two faces is an \emph{edge}. Thus, the boundary of a face consists of edges.

On the positive side, it turns out that almost all the points in $\Omega$ are products. Let us call a subset of $\Omega$ \emph{admissible} if all its points are products. We show:
\begin{itemize}
  \item The interior of $\Omega$ is admissible.
  \item The faces of $\Omega$, excluding the edges, are admissible.
  \item Some but not all edges are admissible.
\end{itemize}
This parallels the situation for the product law for the standard Sobolev spaces $H^s$ (see Theorem \ref{B:Thm1} in the next section).

Concerning the edges, we do not claim to have the optimal result, however. That is, there may be some points on the edges which are products but which are not included in our positive results.

In order to avoid an unduly lengthy paper, we restrict our attention, for the positive results, to the physical space dimension $n=3$, which is of most interest for applications. The cases $n=2$ and $n=1$ will be published in a separate paper (the 2d case is slightly more involved than the 3d case).

Before proceeding to the list of necessary conditions, we make some preliminary observations, and introduce notation and terminology.

It is important to note that if $\Hproduct{s_0 & s_1 & s_2 \\ b_0 & b_1 & b_2}$ is a product, then so is every permutation of its columns. This becomes obvious once we restate~\eqref{A:20} in the following more symmetric form: By Plancherel's theorem and duality,~\eqref{A:20} is equivalent to the trilinear integral estimate
\begin{equation}\label{A:40}
  \Abs{I}
  \lesssim \norm{F_0}\norm{F_1}\norm{F_2},
\end{equation}
where
\begin{equation}\label{A:42}
  I
  =
  \iiint \frac{F_0(X_0)F_1(X_1)F_2(X_2)
  \, \delta(X_0+X_1+X_2) \, \d X_0 \, \d X_1 \, \d X_2}
  {\angles{\xi_0}^{s_0}\angles{\xi_1}^{s_1}\angles{\xi_2}^{s_2}
  \angles{\abs{\tau_0}-\abs{\xi_0}}^{b_0}
  \angles{\abs{\tau_1}-\abs{\xi_1}}^{b_1}
  \angles{\abs{\tau_2}-\abs{\xi_2}}^{b_2} }
\end{equation}
and $X_j=(\tau_j,\xi_j) \in \R^{1+n}$ for $j=0,1,2$. Here $\delta$ is the point mass at $0$ in $\R^{1+n}$, and $\norm{\cdot}$ denotes the $L^2$ norm on $\R^{1+n}$. Without loss of generality we may assume that $F_j \ge 0$ for $j=0,1,2$, hence $I \ge 0$.

Since $\xi_0+\xi_1+\xi_2=0$ in $I$, the triangle inequality implies $\angles{\xi_j} \lesssim \angles{\xi_k} + \angles{\xi_l}$ for all permutations $(j,k,l)$ of $(0,1,2)$, hence the two largest of $\angles{\xi_0}$, $\angles{\xi_1}$ and $\angles{\xi_2}$ are comparable, so we can split
\begin{equation}\label{A:44}
  I = I_{\text{LHH}} + I_{\text{HLH}} + I_{\text{HHL}},
\end{equation}
where the terms on the right hand side are defined by inserting the characteristic functions of the following conditions, respectively, in the integral $I$:
\begin{subequations}\label{A:46}
\begin{align}
  \label{A:46a}
  &\angles{\xi_0} \lesssim \angles{\xi_1} \sim \angles{\xi_2}& &\text{(LHH)}
  \\ 
  \label{A:46b}
  &\angles{\xi_1} \lesssim \angles{\xi_0} \sim \angles{\xi_2}& &\text{(HLH)}
  \\ 
  \label{A:46c}
  &\angles{\xi_2} \lesssim \angles{\xi_0} \sim \angles{\xi_1}& &\text{(HHL)}.
\end{align}
\end{subequations}
Here the mnemonics in the right hand column refer to the relative sizes of the spatial frequencies in the order $(\xi_0,\xi_1,\xi_2)$, with ``L'' and ``H'' standing for low and high frequencies, respectively.

In some situations we also split the $I$'s depending on the signs $\pm_1$ and $\pm_2$ of the temporal frequencies $\tau_1$ and $\tau_2$. Thus,
\begin{equation}\label{A:50}
  I = I^{(+,+)} + I^{(+,-)} + I^{(-,+)} + I^{(-,-)},
\end{equation}
where
\begin{align*}
  &I^{(\pm_1,\pm_2)}
  \\
  &=
  \iiint\limits_{\pm_1\tau_1 \ge 0, \;\pm_2\tau_2 \ge 0}
  \frac{F_0(X_0)F_1(X_1)F_2(X_2)
  \, \delta(X_0+X_1+X_2) \, \d X_0 \, \d X_1 \, \d X_2}
  {\angles{\xi_0}^{s_0}\angles{\xi_1}^{s_1}\angles{\xi_2}^{s_2}
  \angles{\abs{\tau_0}-\abs{\xi_0}}^{b_0}
  \angles{\abs{\tau_1}-\abs{\xi_1}}^{b_1}
  \angles{\abs{\tau_2}-\abs{\xi_2}}^{b_2} }
  \\
  &=
  \iiint\limits_{\pm_1\tau_1 \ge 0, \;\pm_2\tau_2 \ge 0}
  \frac{F_0(X_0)F_1(X_1)F_2(X_2)
  \, \delta(X_0+X_1+X_2) \, \d X_0 \, \d X_1 \, \d X_2}
  {\angles{\xi_0}^{s_0}\angles{\xi_1}^{s_1}\angles{\xi_2}^{s_2}
  \angles{\abs{\tau_0}-\abs{\xi_0}}^{b_0}
  \angles{-\tau_1\pm_1\abs{\xi_1}}^{b_1}
  \angles{-\tau_2\pm_2\abs{\xi_2}}^{b_2} },
\end{align*}
and similarly for $I_{\text{LHH}}$ etc.

In conjunction with the splittings~\eqref{A:44} and~\eqref{A:50}, as well as their combination, it is convenient to use the following rather obvious modifications of the terminology introduced in Definition \ref{A:Def2}: When we say, for instance, that $\left.\Hproduct{s_0 & s_1 & s_2 \\ b_0 & b_1 & b_2}\right\rvert_{\text{HLH}}$ is a product, we mean that~\eqref{A:40} holds for $I_{{\text{HLH}}}$, and if we say that $\left.\Hproduct{s_0 & s_1 & s_2 \\ b_0 & b_1 & b_2}\right\rvert^{(+,+)}_{\text{LHH}}$ is a product, we mean that~\eqref{A:40} holds for $I^{(+,+)}_{\text{LHH}}$, and so on.

We use $x \lesssim y$ as a convenient shorthand for $x \le Cy$, where $C \gg 1$ is a constant which may depend on quantities that are considered fixed. Moreover, $x \sim y$ stands for $x \lesssim y \lesssim x$.

\section{The product law}\label{B}

\subsection{Necessary conditions}

A number of explicit examples given in \S\ref{F} show that any product $\Hproduct{s_0 & s_1 & s_2 \\ b_0 & b_1 & b_2}$ must necessarily satisfy the following $21$ conditions:
\begin{align}
  \label{B:1}
  & b_0 + b_1 + b_2 \ge \frac12
  \\
  \label{B:2}
  & b_0 + b_1 \ge 0
  \\
  \label{B:3}
  & b_0 + b_2 \ge 0
  \\
  \label{B:4}
  & b_1 + b_2 \ge 0
  \\
  \label{B:5}
  & s_0 + s_1 + s_2  \ge \frac{n+1}2 - (b_0 + b_1 + b_2)
  \\
  \label{B:6}
  & s_0 + s_1 + s_2 \ge \frac{n}2 - (b_0 + b_1)
  \\
  \label{B:7}
  & s_0 + s_1 + s_2 \ge \frac{n}2 - (b_0 + b_2)
  \\
  \label{B:8}
  & s_0 + s_1 + s_2 \ge \frac{n}2 - (b_1 + b_2)
  \\
  \label{B:9}
  & s_0 + s_1 + s_2 \ge \frac{n-1}2 - b_0
  \\
  \label{B:10}
  & s_0 + s_1 + s_2 \ge \frac{n-1}2 - b_1
  \\
  \label{B:11}
  & s_0 + s_1 + s_2 \ge \frac{n-1}2 - b_2
  \\
  \label{B:12}
  & s_0 + s_1 + s_2 \ge \frac{n+1}4
  \\
  \label{B:13}
  & (s_0 + b_0) + 2s_1+ 2s_2 \ge \frac{n}2& &\Hproduct{\text{L} & \text{H} & \text{H} \\ & + & -}
  \\
  \label{B:14}
  & 2s_0 + (s_1 + b_1) + 2s_2 \ge \frac{n}2& &\Hproduct{\text{H} & \text{L} & \text{H} \\ + & & -}
  \\
  \label{B:15}
  & 2s_0 + 2s_1 + (s_2 + b_2) \ge \frac{n}2& &\Hproduct{\text{H} & \text{H} & \text{L} \\ + & - &}
  \\
  \label{B:16}
  & s_1 + s_2 \ge -b_0& &\Hproduct{\text{L} & \text{H} & \text{H} \\ & + & +}
  \\
  \label{B:17}
  & s_0 + s_2 \ge -b_1 & &\Hproduct{\text{H} & \text{L} & \text{H} \\ + & & +}
  \\
  \label{B:18}
  & s_0 + s_1 \ge -b_2& &\Hproduct{\text{H} & \text{H} & \text{L} \\ + & + &}
  \\
  \label{B:19}
  & s_1 + s_2 \ge 0& &\Hproduct{\text{L} & \text{H} & \text{H} \\  & \phantom{+}&\phantom{+} }
  \\
  \label{B:20}
  & s_0 + s_2 \ge 0& &\Hproduct{\text{H} & \text{L} & \text{H} \\ \phantom{+} & &\phantom{+} }
  \\
  \label{B:21}
  & s_0 + s_1 \ge 0& &\Hproduct{\text{H} & \text{H} & \text{L} \\ \phantom{+} & \phantom{+}& }.
\end{align}

The tags in the right hand column have the following meaning: The upper row indicates the spatial frequency interaction (LHH, HLH or HHL) in which the condition is necessary. The lower row, if not empty, indicates whether the signs of the respective temporal frequencies are equal (indicated by $+$$+$) or opposite (indicated by $+$$-$). For example,~\eqref{B:13} [resp.~\eqref{B:16}] is needed in the LHH interaction with opposite [resp.~equal] signs for $\tau_1$ and $\tau_2$. An empty lower row means, of course, that the condition is needed regardless of the signs.

The same qualifications are understood to apply also in Theorem \ref{B:Thm2} below.

\begin{definition} Let $\Omega$ be the convex polyhedron in $\R^6$ determined by the above conditions.
\end{definition}

The interior of $\Omega$ corresponds to strict inequality in all the conditions. Each condition determines a face of $\Omega$, corresponding to the case of equality. If at least two of the conditions are equalities, then we are in an edge.

As we said in the introduction, we shall prove (for $n \le 3$) that the interior and the faces minus the edges are admissible. Moreover, some but not all of the edges are admissible.

This parallels the situation for the comparatively trivial product law for the standard Sobolev spaces, which we now recall.

\subsection{Comparison with the product law for $H^s$} This reads as follows:

\begin{theorem}\label{B:Thm1} Let $s_0,s_1,s_2 \in \R$. The product estimate
$$
  \norm{fg}_{H^{-s_0}}
  \le C \norm{f}_{H^{s_1}} \norm{g}_{H^{s_2}}
$$
holds if and only if
\begin{align}
  \label{B:40}
  & s_0 + s_1 + s_2 \ge \frac{n}2
  \\
  \label{B:41}
  & s_0 + s_1 \ge 0
  \\
  \label{B:42}
  & s_0 + s_2 \ge 0
  \\
  \label{B:43}
  & s_1 + s_2 \ge 0
  \\
  \label{B:44}
  &\parbox[t]{10.8 cm}{ If~\eqref{B:40} is an equality, then~\eqref{B:41}--\eqref{B:43} must be strict.}
\end{align}
\end{theorem}

The simple proof of the positive part will be shown later, since the same argument comes up also in the proof of the wave-Sobolev product law. The negative part of the above theorem follows by a standard example which we do not repeat here.

The conditions~\eqref{B:40}--\eqref{B:43} determine a convex polyhedron of points $(s_0,s_1,s_2)$ in $\R^3$. The edges corresponding to equality in~\eqref{B:40} and one of~\eqref{B:41}--\eqref{B:43} are not admissible. On the other hand, the edges corresponding to equality in two of~\eqref{B:41}--\eqref{B:43} are admissible, as long as we stay away from the face given by equality in~\eqref{B:40}. It therefore seems difficult to write down a simple rule telling us which edges are admissible. 

But if instead of talking about edges we talk about equalities, then we can make a simple rule as follows: Replace~\eqref{B:40}--\eqref{B:44} by
\begin{align}
  \label{B:45}
  & s_0 + s_1 + s_2 \ge \frac{n}2
  \\
  \label{B:46}
  & s_0 + s_1 + s_2 \ge \max(s_0,s_1,s_2).
\end{align}
where we combined~\eqref{B:41}--\eqref{B:43} into a single condition. Then~\eqref{B:44} is replaced by the statement:
\begin{equation}\label{B:47}
  \parbox[t]{10.8 cm}{We do not allow both~\eqref{B:45} and~\eqref{B:46} to be equalities.}
\end{equation}

By comparing~\eqref{B:45} and~\eqref{B:46}, both with equality assumed, the last rule can also be reformulated as a list of explicit exceptions as follows:
\begin{align}
  \label{B:48a}
  \parbox[t]{10.8 cm}{If $s_0=\frac{n}2$, then~\eqref{B:40}=\eqref{B:43} must be strict.}
  \\
  \label{B:48b}
  \parbox[t]{10.8 cm}{If $s_1=\frac{n}2$, then~\eqref{B:40}=\eqref{B:42} must be strict.}
  \\
  \label{B:48c}
  \parbox[t]{10.8 cm}{If $s_2=\frac{n}2$, then~\eqref{B:40}=\eqref{B:41} must be strict.}
\end{align}
Here the notation ``\eqref{B:40}=\eqref{B:43}'' indicates that the two conditions coincide.

These ideas help to systematize the much more complicated exceptions along the edges of $\Omega$, which we now discuss.

\subsection{Exceptions on the boundary of $\Omega$}\label{B:50} First rewrite~\eqref{B:1}--\eqref{B:4} as
\begin{align}
  \label{B:51}
  & b_0 + b_1 + b_2 \ge \frac12
  \\
  \label{B:52}
  & b_0 + b_1 + b_2 \ge \max(b_0,b_1,b_2)
\end{align}
Then we impose the rule:
\begin{equation}\label{B:53}
  \parbox[t]{10.8 cm}{We do not allow both~\eqref{B:51} and~\eqref{B:52} to be equalities.}
\end{equation}

Next, consider~\eqref{B:5}--\eqref{B:21}. By symmetry, it suffices to consider the LHH case, hence we ignore those conditions among~\eqref{B:13}--\eqref{B:21} which are not tagged LHH. Moreover, we do not want to compare~\eqref{B:13} with~\eqref{B:16} since they have different sign assumptions, hence we split into $\Hproduct{\text{L} & \text{H} & \text{H} \\ & + & +}$ and $\Hproduct{\text{L} & \text{H} & \text{H} \\ & + & -}$

In the case $\Hproduct{\text{L} & \text{H} & \text{H} \\ & + & +}$ we rewrite the relevant conditions from~\eqref{B:5}--\eqref{B:21} as:
\begin{align}
  \label{B:54}
  & s_0 + s_1 + s_2 \ge \frac{n+1}2 - (b_0 + b_1 + b_2)
  \\
  \label{B:55}
  & s_0 + s_1 + s_2 \ge \frac{n}2 + \max( - b_0 - b_1, - b_0 - b_2, - b_1 - b_2 )
  \\
  \label{B:56}
  & s_0 + s_1 + s_2 \ge \frac{n-1}2 + \max\left( - b_0 , - b_1, - b_2, - \frac{n-3}4 \right)
  \\
  \label{B:57}
  & s_0 + s_1 + s_2 \ge s_0 + \max(0,-b_0)& &\Hproduct{\text{L} & \text{H} & \text{H} \\ & + & +}.
\end{align}
Then we impose the rule:
\begin{equation}\label{B:58}
  \parbox[t]{10.8 cm}{We allow at most one of~\eqref{B:54}--\eqref{B:57} to be an equality.}
\end{equation}

In the case $\Hproduct{\text{L} & \text{H} & \text{H} \\ & + & -}$ we rewrite the relevant conditions from~\eqref{B:5}--\eqref{B:21} as:
\begin{align}
  \label{B:60}
  & s_0 + s_1 + s_2 \ge \frac{n+1}2 - (b_0 + b_1 + b_2)
  \\
  \label{B:61}
  & s_0 + s_1 + s_2 \ge \frac{n}2 + \max( - b_0 - b_1, - b_0 - b_2, - b_1 - b_2 )
  \\
  \label{B:62}
  & s_0 + s_1 + s_2 \ge \frac{n-1}2 + \max\left( - b_0 , - b_1, - b_2, - \frac{n-3}4 \right)
  \\
  \label{B:63}
  & s_0 + s_1 + s_2 \ge \frac{n}4 + \frac{s_0-b_0}2& &\Hproduct{\text{L} & \text{H} & \text{H} \\ & + & -}
  \\
  \label{B:64}
  & s_0 + s_1 + s_2 \ge s_0& &\Hproduct{\text{L} & \text{H} & \text{H} \\ & + & -},
\end{align}
and we impose the rule:
\begin{equation}\label{B:65}
  \parbox[t]{10.8 cm}{We allow at most one of~\eqref{B:60}--\eqref{B:64} to be an equality.}
\end{equation}

An alternative formulation of the above rules is given in Theorem \ref{B:Thm3} below.

The analogous rules for the HLH and HHL cases are obtained by changing the subscript $0$ in the right hand side of~\eqref{B:57},~\eqref{B:63} and~\eqref{B:64} to a $1$ or $2$, respectively.

\subsection{The product law for $H^{s,b}$} We can now formulate the main result:

\begin{theorem}\label{B:Thm2} Let $n = 3$. Assume that $s_0,s_1,s_2,b_0,b_1,b_2 \in \R$ satisfy the conditions~\eqref{B:1}--\eqref{B:21}. Moreover, assume that the rules set out in \S\ref{B:50} are satisfied. Then $\Hproduct{s_0 & s_1 & s_2 \\ b_0 & b_1 & b_2}$ is a product.
\end{theorem}

\begin{remark}
For $n=1$ and $n=2$ the same result holds; the proofs will appear in a separate paper. We expect the same result to hold also for $n \ge 4$. 
\end{remark}

\begin{remark}
In the course of the proof, we break Theorem \ref{B:Thm2} down according to the classification into product types introduced in \S\ref{BB} below, and we restate the theorem in a more explicit form in each case. For practical use, the reader may find these restatements easier to deal with than the general statement in Theorem \ref{B:Thm2}. See \S\S\ref{P:100}--\ref{P:500:0}.
\end{remark}

\begin{remark}
We are not claiming that the boundary rules are necessary, only that they are sufficient. We do expect, however, that~\eqref{B:53} is necessary. This is certainly true in the 1d case, where it can be seen from the standard counterexample for the $H^s$ product law. We also expect~\eqref{B:58} and~\eqref{B:65} to be necessary if all the $b$'s are nonnegative, but if one of the $b$'s is negative, then they can under certain conditions be relaxed somewhat (see Theorem \ref{P:Thm5} below).
\end{remark}

By comparing equalities pairwise within the groups~\eqref{B:51}--\eqref{B:52},~\eqref{B:54}--\eqref{B:57} and~\eqref{B:60}--\eqref{B:64}, we can restate the rules~\eqref{B:53},~\eqref{B:58} and~\eqref{B:65} as an explicit list of exceptions analogous to the list~\eqref{B:48a}--\eqref{B:48c} for the $H^s$ product law:

\begin{theorem}\label{B:Thm3}
Let $n=3$. Assume that~\eqref{B:1}--\eqref{B:21} are verified. Then the rules~\eqref{B:53},~\eqref{B:58} and~\eqref{B:65} for the LHH interaction are equivalent to the following list of exceptions:
\begin{align}
  \label{B:80}
  \parbox[t]{10.8cm}{If $b_0 = \frac12$, then~\eqref{B:1}=\eqref{B:4},~\eqref{B:5}=\eqref{B:8},~\eqref{B:6}=\eqref{B:10} and~\eqref{B:7}=\eqref{B:11} must all be strict.}
  \\
  \label{B:81}
  \parbox[t]{10.8cm}{If $b_1 = \frac12$, then~\eqref{B:1}=\eqref{B:3},~\eqref{B:5}=\eqref{B:7},~\eqref{B:6}=\eqref{B:9} and~\eqref{B:8}=\eqref{B:11} must all be strict.}
  \\
  \label{B:82}
  \parbox[t]{10.8cm}{If $b_2 = \frac12$, then~\eqref{B:1}=\eqref{B:2},~\eqref{B:5}=\eqref{B:6},~\eqref{B:7}=\eqref{B:9} and~\eqref{B:8}=\eqref{B:10} must all be strict.}
  \\
  \label{B:83}
  \parbox[t]{10.8cm}{If $b_0 + b_1 = 1$, then~\eqref{B:5}=\eqref{B:11} must be strict.}
  \\
  \label{B:84}
  \parbox[t]{10.8cm}{If $b_0 + b_2 = 1$, then~\eqref{B:5}=\eqref{B:10} must be strict.}
\\
  \label{B:85}
  \parbox[t]{10.8cm}{If $b_1 + b_2 = 1$, then~\eqref{B:5}=\eqref{B:9} must be strict.}
  \\
  \label{B:86}
  \parbox[t]{10.8cm}{If $b_0 + b_1 = \frac{n-1}4$, then~\eqref{B:6}=\eqref{B:12} must be strict.}
  \\
  \label{B:87}
  \parbox[t]{10.8cm}{If $b_0 + b_2 = \frac{n-1}4$, then~\eqref{B:7}=\eqref{B:12} must be strict.}
  \\
  \label{B:88}
  \parbox[t]{10.8cm}{If $b_1 + b_2 = \frac{n-1}4$, then~\eqref{B:8}=\eqref{B:12} must be strict.}
  \\
  \label{B:89}
  \parbox[t]{10.8cm}{If $b_0 + b_1 + b_2 = \frac{n+1}4$, then~\eqref{B:5}=\eqref{B:12} must be strict.}
  \\
  \label{B:90}
  \parbox[t]{10.8cm}{If $s_0 - b_0 = \frac{n+2}2 - 2(b_0+b_1+b_2)$, then~\eqref{B:5}=\eqref{B:13} must be strict.}
  \\
  \label{B:91}
  \parbox[t]{10.8cm}{If $s_0 - b_0 = \frac{n}2 - 2(b_0+b_1)$, then~\eqref{B:6}=\eqref{B:13} must be strict.}
  \\
  \label{B:92}
  \parbox[t]{10.8cm}{If $s_0 - b_0 = \frac{n}2 - 2(b_0+b_2)$, then~\eqref{B:7}=\eqref{B:13} must be strict.}
  \\
  \label{B:93}
  \parbox[t]{10.8cm}{If $s_0 - b_0 = \frac{n-2}2 - 2b_0$, then~\eqref{B:9}=\eqref{B:13} must be strict.}
  \\
  \label{B:94}
  \parbox[t]{10.8cm}{If $s_0 - b_0 = \frac12$, then~\eqref{B:12}=\eqref{B:13} must be strict.}
  \\
  \label{B:95}
  \parbox[t]{10.8cm}{If $s_0 - b_0 = \frac{n}2 - 2b_0$, then~\eqref{B:19}=\eqref{B:13} must be strict.}
  \\
  \label{B:96}
  \parbox[t]{10.8cm}{If one of~\eqref{B:5}--\eqref{B:12} is an equality, then~\eqref{B:16} and \eqref{B:19} must be strict.}
\end{align}
\end{theorem}

Here the notation ``\eqref{B:1}=\eqref{B:2}'' indicates that the two conditions coincide.

The analogous exceptions for the HLH and HHL cases are obtained by permuting the subscripts in~\eqref{B:90}--\eqref{B:95}.

\subsection{Classification of products}\label{BB}

By permutation invariance, it suffices to prove the main result for products of the following special types:
{
\renewcommand{\theenumi}{\Roman{enumi}}
\begin{enumerate}
\item\label{Type:I} $b_0,b_1,b_2 \ge 0$. Then by symmetry it suffices to consider the subtypes
  \begin{enumerate}
  \item\label{Type:Ia} $b_0=b_1=0 < b_2$,
  \item\label{Type:Ib} $b_0=0 < b_1,b_2$,
  \item\label{Type:Ic} $0 < b_0,b_1,b_2$,
  \end{enumerate}
  \item\label{Type:II} $b_0 < 0 < b_1,b_2$.
\end{enumerate}
}

\subsection{Outline of paper} In \S\ref{F} the counterexamples which imply the necessary conditions are given. In \S\ref{N} we make a dyadic decomposition of the integral $I$ and recall the dyadic estimates which are the fundamental building blocks in the proof of the product laws. We also recall the simple proof of the $H^s$ product law, since that argument is used repeatedly in later sections. The main result, Theorem \ref{B:Thm2}, is proved in \S\S\ref{P:100}--\ref{P:500:0}, broken into sections according to the classification into types as in \S\ref{BB}. In each section we explicitly restate the theorem, and this may be useful also when applying our results, as an alternative to grappling with the general formulation above. The reformulation of the boundary rules, Theorem \ref{B:Thm3}, is proved in \S\ref{Z}.

\section{Counterexamples}\label{F}

To prove the necessity of~\eqref{B:1}--\eqref{B:21} we will estimate the integral $I$, defined by~\eqref{A:42}, on examples of the form $F_0 = \chi_{-C}$, $F_1 = \chi_{A}$ and $F_2 = \chi_{B}$, where $A,B,C \subset \R^{1 + n}$ depend on a parameter $\lambda \gg 1$ and are chosen so that $A + B \subset C$. Since $X_0+X_1+X_2 = 0$ in $I$, this ensures that
$$
  X_1 \in A, \; X_2 \in B \implies X_0=-(X_1+X_2) \in -C,
$$
and then we estimate the weight in $I$ by a power of $\lambda$:
\begin{equation*}
    \angles{\xi_0}^{s_0}\angles{\xi_1}^{s_1}\angles{\xi_2}^{s_2}
  \angles{\abs{\tau_0}-\abs{\xi_0}}^{b_0}
  \angles{\abs{\tau_1}-\abs{\xi_1}}^{b_1}
  \angles{\abs{\tau_2}-\abs{\xi_2}}^{b_2} 
  \sim
  \lambda^{\delta(s_0,s_1,s_2,b_0,b_1,b_2)},
\end{equation*}
where $\delta$ is some linear combination of the $s$'s and $b$'s. Then $I \sim \lambda^{-\delta} \abs{A}\abs{B}$, while $\norm{F_0} \norm{F_1} \norm{F_2}
  \sim
  \abs{A}^{\frac12} \abs{B}^{\frac12} \abs{C}^{\frac12}$.
The estimate~\eqref{A:40} will then imply the condition $\lambda^{\delta}
  \gtrsim
  \abs{A}^{\frac12}\abs{B}^{\frac12}\abs{C}^{-\frac12}$.
If we have an estimate of the form
\begin{equation} \label{eq:24}
  \frac{\abs{A}^{\frac12} \abs{B}^{\frac12}}{\abs{C}^{\frac12}}
  \sim
  \lambda^{d(n)},
\end{equation}
then we deduce the necessary condition $\delta=\delta(s_0,s_1,s_2,b_0,b_1,b_2) \ge d(n)$.

In the following we split $\xi \in \R^n$ as $\xi=(\xi_1,\xi')$, where $\xi'=(\xi_2,\dots,\xi_n) \in \R^{n-1}$. To avoid any confusion, we emphasize that in this notation the subscript refers to coordinates, whereas elsewhere we use subscripts to label different vectors.

\subsection{Necessity of~\eqref{B:1}} This is obtained by scaling only the temporal variables:
  \begin{align*}
    & A = B = \setb{(\tau, \xi)}{%
      \lambda \le \tau \le 2 \lambda,
      \abs{\xi} \le 1
    }, & \abs{A} = \abs{B} \sim \lambda, \\
    & C = \setb{(\tau, \xi)}{%
      2 \lambda \le \tau \le 4 \lambda,
      \abs{\xi} \le 2
    }, & \abs{C} \sim \lambda, \\
    & \delta = b_0+b_1+b_2, & d = \frac{1}2.
  \end{align*}
  
\subsection{Necessity of~\eqref{B:2}--\eqref{B:4}} By symmetry, it suffices to show~\eqref{B:3}, and for this we choose:
  \begin{align*}
    & A = \setb{(\tau, \xi)}{%
      \abs{\tau} \le 1,
      \abs{\xi} \le 1
    }, & \abs{A} \sim 1, \\
    & B = \setb{(\tau, \xi)}{%
      \abs{\tau - \lambda} \le 1,
      \abs{\xi} \le 1
    }, & \abs{B} \sim 1, \\
    & C = \setb{(\tau, \xi)}{%
      \abs{\tau - \lambda} \le 2,
      \abs{\xi} \le 2
    }, & \abs{C} \sim 1, \\
    & \delta = b_0+b_2, & d = 0.
  \end{align*}

\subsection{Necessity of~\eqref{B:5}} This is obtained by scaling all variables:
  \begin{align*}
    & A = B = \setb{(\tau, \xi)}{%
      \abs{\tau} \le \frac{\lambda}2,
      \lambda \le \xi_1 \le 2 \lambda,
      \abs{\xi'} \le \lambda
    }, & \abs{A} = \abs{B} \sim \lambda^{n + 1}, \\
    & C = \setb{(\tau, \xi)}{%
      \abs{\tau} \le \lambda,
      2 \lambda \le \xi_1 \le 4 \lambda,
      \abs{\xi'} \le 2 \lambda
    }, & \abs{C} \sim \lambda^{n + 1}, \\
    & \delta = s_0+s_1+s_2+b_0+b_1+b_2, & d = \frac{n + 1}2.
  \end{align*}

\subsection{Necessity of~\eqref{B:6}--\eqref{B:8}} By symmetry, it suffices to show~\eqref{B:7}, and for this we choose:
  \begin{align*}
    & A = \setb{(\tau, \xi)}{%
      \abs{\tau - \abs{\xi}} \le 1,
      \lambda \le \xi_1 \le \frac{5\lambda}4,
      \abs{\xi'} \le \frac{\lambda}4
    }, & \abs{A} \sim \lambda^n, \\
    & B = \setb{(\tau, \xi)}{%
      \abs{\tau} \le \frac{\lambda}2,
      2 \lambda \le \xi_1 \le \frac{5\lambda}2,
      \abs{\xi'} \le \frac{\lambda}2
    }, & \abs{B} \sim \lambda^{n + 1}, \\
    & C = \setb{(\tau, \xi)}{%
      \abs{\tau} \le \frac{5\lambda}2,
      3 \lambda \le \xi_1 \le 4 \lambda,
      \abs{\xi'} \le \lambda
    }, & \abs{C} \sim \lambda^{n + 1}, \\
    & \delta = s_0+s_1+s_2+b_0+b_2, & d = \frac{n}{2}.
  \end{align*}
  
\subsection{Necessity of~\eqref{B:9}--\eqref{B:11}} By symmetry, it suffices to show~\eqref{B:9}, and for this we choose:
  \begin{align*}
    & A = \setb{(\tau, \xi)}{%
      \abs{\tau - \abs{\xi}} \le 1,
      \lambda \le \xi_1 \le \frac{3\lambda}2,
      \abs{\xi'} \le \frac{\lambda}2
    }, & \abs{A} \sim \lambda^n, \\
    & B = \setb{(\tau, \xi)}{%
      \abs{\tau + \abs{\xi}} \le 1,
      \lambda \le \xi_1 \le \frac{3\lambda}2,
      \abs{\xi'} \le \frac{\lambda}2
    }, & \abs{B} \sim \lambda^n, \\
    & C = \setb{(\tau, \xi)}{%
      \abs{\tau} \le \frac{3\lambda}2,
      2 \lambda \le \xi_1 \le 3 \lambda,
      \abs{\xi'} \le \lambda
    }, & \abs{C} \sim \lambda^{n + 1}, \\
    & \delta = s_0+s_1+s_2+b_0, & d = \frac{n - 1}{2}.
  \end{align*}

\subsection{Necessity of~\eqref{B:12}} This represents the effect of Lorentz transformations (concentration along null directions):
  \begin{align*}
    & A = B = \setb{(\tau, \xi)}{%
      \abs{\tau - \xi_1} \le 1,
      \lambda \le \xi_1 \le 2 \lambda,
      \abs{\xi'} \le \sqrt{\lambda}
    }, & \abs{A} = \abs{B} \sim \lambda^{\frac{n + 1}2}, \\
    & C = \setb{(\tau, \xi)}{%
      \abs{\tau - \xi_1} \le 2,
      2 \lambda \le \xi_1 \le 4 \lambda,
      \abs{\xi'} \le 2 \sqrt{\lambda}
    }, & \abs{C} \sim \lambda^{\frac{n + 1}2}, \\
    & \delta = s_0+s_1+s_2, & d = \frac{n + 1}4.
  \end{align*}

\subsection{Necessity of~\eqref{B:13}--\eqref{B:15}} By symmetry, it suffices to show~\eqref{B:13}, and for this we choose:
  \begin{align*}
    & A = \setb{(\tau, \xi)}{%
      \abs{\tau - \xi_1} \le 1,
      \abs{\xi_1 - \lambda^2} \le \lambda,
      \abs{\xi'} \le \lambda, 
      \xi_2 \ge \frac12 \lambda
    }, & \abs{A} \sim \lambda^n, \\
    & B = \setb{(\tau, \xi)}{%
      \abs{\tau - \xi_1} \le 1,
      \abs{\xi_1 + \lambda^2} \le \lambda,
      \abs{\xi'} \le \lambda,
      \xi_2 \ge \frac12 \lambda
    }, & \abs{B} \sim \lambda^n, \\
    & C = \setb{(\tau, \xi)}{%
      \abs{\tau - \xi_1} \le 2,
      \abs{\xi_1} \le 2 \lambda,
      \abs{\xi'} \le 2 \lambda,
      \xi_2 \ge \lambda
    }, & \abs{C} \sim \lambda^n, \\
    & \delta = s_0 + b_0 + 2 s_1 + 2 s_2, & d = \frac{n}{2}.
  \end{align*}

\subsection{Necessity of~\eqref{B:16}--\eqref{B:18}} By symmetry, it suffices to show~\eqref{B:16}, and for this we choose:
  \begin{align*}
    & A = \setb{(\tau, \xi)}{%
      \abs{\tau - \lambda} \le 1,
      \abs{\xi_1 - \lambda} \le 1,
      \abs{\xi'} \le 1
    }, & \abs{A} \sim 1, \\
    & B = \setb{(\tau, \xi)}{%
      \abs{\tau - \lambda} \le 1,
      \abs{\xi_1 + \lambda} \le 1,
      \abs{\xi'} \le 1
    }, & \abs{B} \sim 1, \\
    & C = \setb{(\tau, \xi)}{%
      \abs{\tau - 2 \lambda} \le 2,
      \abs{\xi_1} \le 2,
      \abs{\xi'} \le 2
    }, & \abs{C} \sim 1, \\
    & \delta = s_1+s_2+b_0, & d = 0.
  \end{align*}
  
\subsection{Necessity of~\eqref{B:19}--\eqref{B:21}} By symmetry, it suffices to show~\eqref{B:20}, and for this we choose:
  \begin{align*}
    & A = \setb{(\tau, \xi)}{%
      \abs{\tau} \le 1,
      \abs{\xi_1} \le 1,
      \abs{\xi'} \le 1
    }, & \abs{A} \sim 1, \\
    & B = \setb{(\tau, \xi)}{%
      \abs{\tau} \le 1,
      \abs{\xi_1 - \lambda} \le 1,
      \abs{\xi'} \le 1
    }, & \abs{B} \sim 1, \\
    & C = \setb{(\tau, \xi)}{%
      \abs{\tau} \le 2,
      \abs{\xi_1 - \lambda} \le 2,
      \abs{\xi'} \le 2
    }, & \abs{C} \sim 1, \\
    & \delta = s_0+s_2, & d = 0.
  \end{align*}
  
\section{Notation and preliminaries}\label{N}

\subsection{Dyadic decompositions}\label{N:1}

Throughout, $M$, $N$ and $L$, as well as their indexed counterparts, denote dyadic numbers of the form $2^j$, $j \in \{0,1,2,\dots\}$. We rely on dyadic decompositions with respect to the size of the weights in the $H^{s,b}$-norm. In some cases we also decompose with respect to the sign of the temporal frequency. 

Given $u \in H^{s,b}$, we define the $L^2$-function $F \ge 0$ by
\begin{equation}\label{N:2}
  F(X) = \angles{\xi}^s\angles{\abs{\tau}-\abs{\xi}}^b\widetilde u(X),
\end{equation}
where $X = (\tau,\xi)$. We shall use the shorthand
$F^{N}(X) = \chi_{\angles{\xi} \sim N} F(X)$, $F^{N,L}(X) = \chi_{\angles{\abs{\tau}-\abs{\xi}} \sim L} F^N(X)$ and $F^{N,L,\pm}(X) = \chi_{\pm \tau \ge 0} F^{N,L}(X)$, and correspondingly we then define $u^{N}$, $u^{N,L}$ and $u^{N,L,\pm}$ as in \eqref{N:2}, replacing $F$ there by $F^{N}$, $F^{N,L}$ and $F^{N,L,\pm}$, respectively. Note that $\sum_N \norm{F^{N}}^2 \sim \norm{F}^2$, $\sum_L \norm{F^{N,L}}^2 \sim \norm{F^N}^2$ and $\sum_L \norm{F^{N,L,\pm}}^2 \lesssim \norm{F^N}^2$.

Defining the trilinear convolution form
$$
  \triple{F_0}{F_1}{F_2}
  =
  \iiint
  F_0(X_0) F_1(X_1) F_2(X_2)
  \, \delta(X_0+X_1+X_2) \, \d X_0 \, \d X_1 \, \d X_2,
$$
we then obtain
\begin{equation}\label{N:20}
  I
  \lesssim
  \sum_{\Nbold,\Lbold}
  \frac{
  \triple{F_0^{N_0,L_0}}{F_1^{N_1,L_1}}{F_2^{N_2,L_2}}
  }
  {N_0^{s_0}N_1^{s_1}N_2^{s_2} L_0^{b_0}L_1^{b_1}L_2^{b_2}},
\end{equation}
where $I$ is given by \eqref{A:42} and we set $\Nbold=(N_0,N_1,N_2)$ and $\Lbold=(L_0,L_1,L_2)$. We use the shorthand $\Nmin^{012} = \min(N_0,N_1,N_2)$,
and similarly for the $L$'s, and for other indexes than $012$.

We also have the analogues of~\eqref{N:20} for $I_{\text{LHH}}$, $I_{\text{HLH}}$ and $I_{\text{HLH}}$, obtained by inserting the characteristic functions of the following conditions, respectively, in the sum on the right hand side of~\eqref{N:20}: $N_0 \le N_1 \sim N_2$ (LHH), $N_1 \le N_0 \sim N_2$ (HLH) and $N_2 \le N_0 \sim N_1$ (HHL).

Note that if $1 \le A < B$ and $a \in \R$, then
\begin{equation}\label{N:36}
  \sum_{A \le L \le B} L^a
  \sim
  \begin{cases}
  B^a &\text{if $a > 0$}
  \\
  \log\angles{\frac{B}{A}}
  &\text{if $a = 0$}
  \\
  A^a &\text{if $a < 0$}.
  \end{cases}
\end{equation}
We frequently apply the estimate, for any $\varepsilon > 0$,
\begin{equation}\label{N:38}
  \log\angles{B} \le C_\varepsilon B^\varepsilon \qquad \text{for all $B \ge 1$}.
\end{equation}

\subsection{Hyperbolic Leibniz rule}\label{N:90}

We recall a well-known ``Leibniz rule'' for hyperbolic weights (a proof can be found, for example, in \cite[Lemma 3.4]{Selberg:2002b}): Assume that $\tau_0+\tau_1+\tau_2=0$ and $\xi_0+\xi_1+\xi_2=0$, as in the integral $I$, and let $\pm_1$ and $\pm_2$ denote the signs of $\tau_1$ and $\tau_2$, respectively. Then
\begin{equation}\label{N:92}
  \bigabs{\abs{\tau_0}-\abs{\xi_0}}
  \lesssim
  \bigabs{-\tau_1\pm_1\abs{\xi_1}}
  +
  \bigabs{-\tau_2\pm_2\abs{\xi_2}}
  +
  \mathfrak b_{(\pm_1,\pm_2)}(\xi_0,\xi_1,\xi_2),
\end{equation}
where
\begin{equation}\label{N:94}
  \mathfrak b_{(\pm_1,\pm_2)}(\xi_0,\xi_1,\xi_2)
  =
  \begin{cases}
  \abs{\xi_1}+\abs{\xi_2}-\abs{\xi_0} &\text{if $\pm_1=\pm_2$}
  \\
  \abs{\xi_0}-\bigabs{\abs{\xi_1}-\abs{\xi_2}} &\text{if $\pm_1\neq\pm_2$}.
  \end{cases}
\end{equation}
Note the estimate
\begin{equation}\label{N:96}
  \mathfrak b_{(\pm_1,\pm_2)}(\xi_0,\xi_1,\xi_2)
  \lesssim
  \begin{cases}
  \min(\abs{\xi_1},\abs{\xi_2}) &\text{if $\pm_1=\pm_2$}
  \\
  \min(\abs{\xi_0},\abs{\xi_1},\abs{\xi_2}) &\text{if $\pm_1\neq\pm_2$}.
  \end{cases}
\end{equation}
We define associated bilinear operators $\mathfrak B_{(\pm_1,\pm_2)}^\alpha$ by
\begin{multline}\label{N:97}
  \mathcal F\left\{\mathfrak B_{(\pm_1,\pm_2)}^\alpha(f,g)\right\}(\xi_0)
  \\
  =
  \iint
  \left(\mathfrak b_{(\pm_1,\pm_2)}(\xi_0,\xi_1,\xi_2)\right)^\alpha
  \widehat f(\xi_1) \widehat g(\xi_2)
  \,\delta(\xi_0+\xi_1+\xi_2)
  \d\xi_1 \d\xi_2
\end{multline}
for $f,g \in \mathcal S(\R^n)$, where $\mathcal F f = \widehat f$ denotes the Fourier transform.

\subsection{The dyadic building blocks}\label{N:40}

On the one hand, we have the more or less trivial ``Sobolev type'' estimate, which has the same form in all dimensions:
\begin{equation}\label{N:42}
  \triple{F_0^{N_0}}{F_1^{N_1}}{F_2^{N_2,L_2}}
  \lesssim
  \left[ \bigl(\Nmin^{012}\bigr)^n L_2 \right]^{\frac12}
  \bignorm{F_0^{N_0}}
  \bignorm{F_1^{N_1}}
  \bignorm{F_2^{N_2,L_2}}.
\end{equation}
By the Cauchy-Schwarz inequality, this can be reduced to a volume estimate; see, for example, \cite[Eq.\ (37)]{Tao:2001}.

On the other hand, there are the much deeper ``wave type'' estimates which in the case $n=3$ look as follows:
\begin{equation}\label{N:50}
  \triple{F_0^{N_0}}{F_1^{N_1,L_1,\pm_1}}{F_2^{N_2,L_2,\pm_2}}       
  \le C
  \bignorm{F_0^{N_0}}
  \bignorm{F_1^{N_1,L_1,\pm_1}}
  \bignorm{F_2^{N_2,L_2,\pm_2}}
\end{equation}
holds with
\begin{align}
  \label{N:54}
  C
  &\sim \left( \Nmin^{012} \Nmin^{12} L_1 L_2 \right)^{\frac12}
  \\
  \label{N:56}
  C    
  &\sim \left( N_0^2 L_1 L_2 \right)^{\frac12} \qquad \qquad \text{if $N_0 \ll N_1 \sim N_2$ and $\pm_1=\pm_2$}.
\end{align}
This follows via a transfer principle (see, e.g., Lemma 4 in \cite{Selberg:2007d}) from analogous estimates for the homogeneous wave equation (see Theorem 12.1 in~\cite{Foschi:2000}; estimates of this type were first investigated in \cite{Klainerman:1993, Klainerman:1996}).

\subsection{A summation argument and the proof of the $H^s$ product law}\label{N:98}

As a warm-up for the proof of the main result, we prove the $H^s$ product law. By a dyadic decomposition, and using the notation above but temporarily reducing $X \in \R^{1+n}$ to just $\xi \in \R^n$, we need to prove
\begin{equation}\label{N:100}
  \sum_{\Nbold}
  \frac{\triple{F_0^{N_0}}{F_1^{N_1}}{F_2^{N_2}}}
  {N_0^{s_0}N_1^{s_1}N_2^{s_2}}
  \lesssim
  \norm{F_0}
  \norm{F_1}
  \norm{F_2},
\end{equation}
where for the moment $\norm{\cdot}$ is the norm on $L^2(\R^n)$ instead of $L^2(\R^{1+n})$. Now we apply the ``Sobolev type'' estimate
\begin{equation}\label{N:102}
  \triple{F_0^{N_0}}{F_1^{N_1}}{F_2^{N_2}}
  \lesssim
  \bigl(\Nmin^{012}\bigr)^{\frac{n}2}
  \bignorm{F_0^{N_0}}
  \norm{F_1^{N_2}}
  \bignorm{F_2^{N_2}},
\end{equation}
whose proof essentially reduces (see \cite{Tao:2001}) to the fact that each $F_k^{N_k}$ is supported in a ball in $\R^n$ of radius comparable to $N_k$. By symmetry it suffices to consider the LHH interaction, so we are left with the sum
$$
  S 
  =
  \sum_{\Nbold} \chi_{N_0 \le N_1 \sim N_2}
  \frac{ N_0^{\frac{n}2} }{ N_0^{s_0} N_1^{s_1} N_2^{s_2} }
  \bignorm{F_0^{N_0}}
  \bignorm{F_1^{N_1}}
  \bignorm{F_2^{N_2}}.
$$
Setting $A=\frac{n}2-s_0$ and $B=s_1+s_2$ for the sake of generality (we shall reuse the following argument several times), we then have
\begin{equation}\label{D:40}
\begin{aligned}
  S
  &\lesssim
  \sum_{\Nbold}
  \chi_{N_0 \le N_1 \sim N_2}
  \frac{ N_0^{A} }{ N_1^{B} }
  \bignorm{F_0^{N_0}}
  \bignorm{F_1^{N_1}}
  \bignorm{F_2^{N_2}}
  \\
  &\lesssim
  \norm{F_0}
  \sum_{N_1,N_2} \chi_{N_1 \sim N_2}
  \frac{\Sigma_A(N_1)}{N_1^{B}}
  \bignorm{F_1^{N_1}}
  \bignorm{F_2^{N_2}},
\end{aligned}
\end{equation}
where
\begin{equation}\label{D:42}
  \Sigma_A(N_1) =
  \sum_{N_0} \chi_{N_0 \le N_1} N_0^{A}
  \sim
  \begin{cases}
  N_1^{A} &\text{if $A > 0$},
  \\
  \log\angles{N_1} &\text{if $A = 0$},
  \\
  1 &\text{if $A < 0$}.
  \end{cases}
\end{equation}
The estimate $S \lesssim \norm{F_0}\norm{F_1}\norm{F_2}$ now follows by the Cauchy-Schwarz inequality provided that (i) $B \ge A$, (ii) $B \ge 0$ and (iii) we exclude $A=B=0$, as these conditions guarantee that $\Sigma_A(N_1) \lesssim N_1^B$.

This proves, in particular, the positive part of Theorem~\ref{B:Thm1}. The same argument will be reused many times in the proof of the main result, which we now begin.

For the remainder of the paper we assume $n=3$.

\section{The case $b_0=b_1=0 < b_2$}\label{P:100}

Then the product law simplifies to:

\begin{theorem}\label{P:Thm1} Let $n = 3$. Set $b_0=b_1=0$ and assume that
\begin{align}
  \label{P:101}
  &b_2 > \frac12
  \\
  \label{P:102}
  &s_0 + s_1 + s_2 \ge \frac32
  \\
  \label{P:103}
  &s_0 + s_1 + s_2 \ge \max(s_0,s_1,s_2),
\end{align}
and that ~\eqref{P:102} and~\eqref{P:103} are not both equalities. Then $\Hproduct{s_0 & s_1 & s_2 \\ 0 & 0 & b_2}$ is a product.
\end{theorem}

By dyadic decomposition we reduce to proving
\begin{equation}\label{P:110}
  \sum_{\Nbold}
  \frac{S_{\Nbold}}
  {N_0^{s_0}N_1^{s_1}N_2^{s_2}}
  \lesssim \norm{F_0} \norm{F_1} \norm{F_2},
\end{equation}
where $S_{\Nbold}
  =
  \sum_{L_2}
  L_2^{-b_2}
  \triple{F_0^{N_0}}{F_1^{N_1}}{F_2^{N_2,L_2}}$. In fact, it is enough to prove
\begin{equation}\label{P:114}
  S_{\Nbold}
  \lesssim
  \bigl(\Nmin^{012}\bigr)^{\frac32}
  \bignorm{F_0^{N_0}}
  \bignorm{F_1^{N_1}}
  \bignorm{F_2^{N_2}},
\end{equation}
since then~\eqref{P:110} follows by the argument used to prove the $H^s$ product law in \S\ref{N:98}. But \eqref{N:42} implies \eqref{P:114} up to multiplication of the right hand side by $\sum_{L_2} L_2^{\frac12-b_2}$, which converges since $b_2 > \frac12$.

\section{The case $b_0=0 < b_1,b_2$}\label{P:200}

Then the product law reads:

\begin{theorem}\label{P:Thm2} Let $n = 3$. Set $b_0=0$ and assume
\begin{align}
  \label{P:201}
  &b_1, b_2 > 0
  \\
  \label{P:202}
  &b_1 + b_2 \ge \frac12
  \\
  \label{P:203}
  &s_0 + s_1 + s_2 \ge 2 - (b_1 + b_2)
  \\
  \label{P:204}
  &s_0 + s_1 + s_2 \ge \frac32 - b_1
  \\
  \label{P:205}
  &s_0 + s_1 + s_2 \ge \frac32 - b_2
  \\
  \label{P:206}
  &s_0 + s_1 + s_2 \ge 1
  \\
  \label{P:207}
  &s_0 + 2(s_1 + s_2) \ge \frac32& &\Hproduct{\text{L} & \text{H} & \text{H} \\ & + & -}
  \\
  \label{P:208}
  &s_1 + s_2 \ge 0& &\Hproduct{\text{L} & \text{H} & \text{H} \\ & \phantom{+} & \phantom{+}}
  \\
  \label{P:209}
  &s_0 + s_2 \ge 0& &\Hproduct{\text{H} & \text{L} & \text{H} \\ \phantom{+} & &\phantom{+}}
  \\
  \label{P:210}
  &s_0 + s_1 \ge 0& &\Hproduct{\text{H} & \text{H} & \text{L} \\ \phantom{+} & \phantom{+} &},
\end{align}
as well as the exceptions:
\begin{align}
  \label{P:211a}
  &\parbox[t]{10.8cm}{If $b_1=\frac12$, then~\eqref{P:203}=\eqref{P:205} must be strict.}
    \\
  \label{P:211b}
  &\parbox[t]{10.8cm}{If $b_1=\frac12$, then~\eqref{P:204}=\eqref{P:206} must be strict.}
  \\
  \label{P:212a}
  &\parbox[t]{10.8cm}{If $b_2=\frac12$, then~\eqref{P:203}=\eqref{P:204} must be strict.}
  \\
  \label{P:212b}
  &\parbox[t]{10.8cm}{If $b_2=\frac12$, then~\eqref{P:205}=\eqref{P:206} must be strict.}
  \\
  \label{P:215}
  &\parbox[t]{10.8cm}{If $b_1+b_2=1$, then~\eqref{P:203}=\eqref{P:206} must be strict.}
  \\
  \label{P:216}
  &\parbox[t]{10.8cm}{We require~\eqref{P:207} to be strict if $s_0$ takes one of the values $\frac12$, $\frac32$, $\frac32 - 2b_1$, $\frac32 - 2b_2$ or $\frac52 - 2(b_1+b_2)$.}
  \\ 
  \label{P:217}
  &\parbox[t]{10.8cm}{If one of~\eqref{P:203}--\eqref{P:206} is an equality, then~\eqref{P:208}--\eqref{P:210} must be strict.}
\end{align}
Then $\Hproduct{s_0 & s_1 & s_2 \\ 0 & b_1 & b_2}$ is a product.
\end{theorem}

By dyadic decomposition we reduce to proving \eqref{P:110} for
\begin{equation}\label{P:222}
  S_{\Nbold}
  =
  \sum_{\Lbold} \chi_{L_1 \le L_2}
  \frac{
  \triple{F_0^{N_0}}{F_1^{N_1,L_1}}{F_2^{N_2,L_2}}
  }
  {L_1^{b_1} L_2^{b_2}},
\end{equation}
where $\Lbold=(L_1,L_2)$. Here the assumption $L_1 \le L_2$ is justified by symmetry.

For the moment we shall assume strict inequality in \eqref{P:202}:
\begin{equation}\label{P:224}
  b_1+b_2 > \frac12.
\end{equation}
The case $b_1+b_2=\frac12$ is handled later, in \S\ref{P:280}.

\subsection{The HLH case}\label{P:230}

Here we assume $N_1 \le N_0 \sim N_2$. By~\eqref{N:42}--\eqref{N:54} we then know that \eqref{N:50} holds with
\begin{equation}\label{P:231}
  C \sim N_1 L_1^{\frac12} \left[ \min( N_1, L_2 ) \right]^{\frac12}.
\end{equation}
To resolve the minimum, we split into the subcases $L_2 \le N_1$ and $L_2 > N_1$.

\subsubsection{The subcase $L_2 \le N_1$}\label{P:233} Applying this restriction in \eqref{P:222}, we get
\begin{equation}\label{P:234}
  S_{\Nbold}
  \lesssim
  N_1\sigma_{\frac12}(N_1)
  \bignorm{F_0^{N_0}}
  \bignorm{F_1^{N_1}}
  \bignorm{F_2^{N_2}},
\end{equation}
where we write, for $p > 0$,
\begin{equation}\label{P:236}
  \sigma_p(N)
  =
  \sum_{\Lbold}
  \chi_{L_1 \le L_2 \le N} L_1^{\frac12-b_1} L_2^{p-b_2}.
\end{equation}
Using \eqref{N:36} (with $A=1$) repeatedly we find
\begin{equation}\label{P:238}
\begin{aligned}
  \sigma_p(N)
  &\sim
  \sum_{L_2 \le N}
  L_2^{p-b_2}
  \times
  \begin{cases}
  L_2^{\frac12-b_1} &\text{if $b_1 < \frac12$}
  \\
  \log\angles{L_2} &\text{if $b_1 = \frac12$}
  \\
  1 &\text{if $b_1 > \frac12$}
  \end{cases}
  \\
  &\lesssim
  \begin{cases}
  N^{\frac12+p-b_1-b_2} &\text{if $b_1 < \frac12$, $b_1+b_2 < \frac12+p$}
  \\
  \log\angles{N} &\text{if $b_1 < \frac12$, $b_1+b_2 = \frac12+p$}
  \\
  1 &\text{if $b_1 < \frac12$, $b_1+b_2 > \frac12+p$}
  \\
  N^{p-b_2}\log\angles{N} &\text{if $b_1 = \frac12$, $b_2 < p$}
  \\
  \log^2\angles{N} &\text{if $b_1 = \frac12$, $b_2 = p$}
  \\
  1 &\text{if $b_1 = \frac12$, $b_2 > p$}
  \\
  N^{p-b_2} &\text{if $b_1 > \frac12$, $b_2 < p$}
  \\
  \log\angles{N} &\text{if $b_1 > \frac12$, $b_2 = p$}
  \\
  1 &\text{if $b_1 > \frac12$, $b_2 > p$}.
  \end{cases}
\end{aligned}
\end{equation}
Applying this to \eqref{P:234} yields
\begin{equation}\label{P:240}
  \sum_{\Nbold}
  \chi_{N_1 \le N_2 \sim N_0}
  \frac{S_{\Nbold}}
  {N_0^{s_0}N_1^{s_1}N_2^{s_2}}
  \lesssim
  \norm{F_1}
  \sum_{\Nbold}
  \chi_{N_1 \le N_2 \sim N_0}
  \frac{N_1^A}{N_0^B}
  \bignorm{F_0^{N_0}}
  \bignorm{F_2^{N_2}},
\end{equation}
where $A$ depends on the $b$'s, whereas $B=s_0+s_2$ in all cases. Therefore, by the argument shown in \S\ref{N:98}, it suffices to check that $B \ge A$ and $B \ge 0$, and moreover that $A=B=0$ cannot happen. Note that $B \ge 0$ is the same as \eqref{P:209}. Logarithmic factors are estimated as in \eqref{N:38}.

\begin{itemize}
\item
If $\sigma_{\frac12}(N_1) \lesssim (N_1)^{1-b_1-b_2}$, then $A=2-s_1-b_1-b_2$, hence $B \ge A$ is~\eqref{P:203}, and~\eqref{P:217} excludes $A=B=0$.
\item
If $\sigma_{\frac12}(N_1) \lesssim (N_1)^{\frac12-b_2}$, then $A=\frac32-s_1-b_2$, hence $B \ge A$ is~\eqref{P:205}, and~\eqref{P:217} excludes $A=B=0$.
\item
If $\sigma_{\frac12}(N_1) \lesssim 1$, then $A=1-s_1$, so $B \ge A$ is~\eqref{P:206}, and~\eqref{P:217} excludes $A=B=0$.
\item If $\sigma_{\frac12}(N_1) \lesssim (N_1)^{\varepsilon}$, then $A=1-s_1+\varepsilon$. But now either $b_2=\frac12$ or $b_1+b_2=1$,  so~\eqref{P:212b} or~\eqref{P:215}, respectively, guarantee that~\eqref{P:206} is strict, hence $B > A$ for $\varepsilon > 0$ small enough.
\item
If $\sigma_{\frac12}(N_1) \lesssim (N_1)^{\frac12-b_2+\varepsilon}$, then $A=\frac32-s_1-b_2+\varepsilon$. But now $b_1 = \frac12$, so~\eqref{P:211a} implies that \eqref{P:205} is strict, hence $B > A$ for $\varepsilon > 0$ small enough.
\end{itemize}

\subsubsection{The subcase $L_2 > N_1$}\label{P:241} Restricting \eqref{P:222} accordingly, and noting that \eqref{N:50} now holds with $C^2 \sim N_1^3 L_1$, we get
\begin{equation}\label{P:242}
  S_{\Nbold}
  \lesssim
  N_1^{\frac32}\gamma(N_1)
  \bignorm{F_0^{N_0}}
  \bignorm{F_1^{N_1}}
  \bignorm{F_2^{N_2}},
\end{equation}
where
\begin{equation}\label{P:244}
  \gamma(N)
  =
  \sum_{\Lbold}
  \chi_{L_1 \le L_2} \chi_{L_2 \ge N} L_1^{\frac12-b_1} L_2^{-b_2}.
\end{equation}
Recalling that $b_1+b_2 > \frac12$ and $b_2 > 0$, by~\eqref{P:224} and~\eqref{P:201}, we find
\begin{equation}\label{P:246}
\begin{aligned}
  \gamma(N)
  &\sim
  \sum_{L_2 \ge N}
  L_2^{-b_2}
  \times
  \begin{cases}
  L_2^{\frac12-b_1} &\text{if $b_1 < \frac12$}
  \\
  \log\angles{L_2} &\text{if $b_1 = \frac12$}
  \\
  1 &\text{if $b_1 > \frac12$}
  \end{cases}
  \\
  &\lesssim
  \begin{cases}
  N^{\frac12-b_1-b_2} &\text{if $b_1 < \frac12$}
  \\
  N^{\varepsilon-b_2} &\text{if $b_1 = \frac12$, and for any $\varepsilon > 0$}
  \\
  N^{-b_2} &\text{if $b_1 > \frac12$}.
  \end{cases}
\end{aligned}
\end{equation}
Applying this in \eqref{P:242} we again get \eqref{P:240} with $B=s_0+s_2$, and with the choice of $A$ depending on the $b$'s. So it is enough to check that $B \ge A$ and $B \ge 0$, and moreover that $A=B=0$ cannot happen. Note that $B \ge 0$ is the same as \eqref{P:209}.

\begin{itemize}
\item
If $b_1 < \frac12$, then $A=2-s_1-b_1-b_2$, hence $B \ge A$ is~\eqref{P:203}, and~\eqref{P:217} excludes $A=B=0$.
\item
If $b_1 = \frac12$, then $A=\frac32-s_1-b_2+\varepsilon$ for any $\varepsilon > 0$. But~\eqref{P:211a} implies \eqref{P:205} strict, hence $B > A$ for $\varepsilon > 0$ small enough.
\item
If $b_1 > \frac12$, then $A=\frac32-s_1-b_2$, hence $B \ge A$ is~\eqref{P:205}, and~\eqref{P:217} excludes $A=B=0$.
\end{itemize}

\begin{remark}\label{P:Remark1}
Note that the conditions~\eqref{P:207},~\eqref{P:208} and ~\eqref{P:210} were not used in the HLH case, and moreover we did not use~\eqref{P:204}, due to the symmetry assumption $L_1 \le L_2$. These conditions can therefore also be deleted from the exceptional rules, and in particular~\eqref{P:211b},~\eqref{P:212a} and \eqref{P:216} are not needed at all. 
\end{remark}

\subsection{The HHL case} This works by an obvious modification of the argument for the HLH case, and the analogue Remark \ref{P:Remark1} remains valid (now it is the conditions~\eqref{P:207},~\eqref{P:209},~\eqref{P:210} and \eqref{P:204} that are not needed, with the corresponding changes to the exceptional rules).

\subsection{The HHL case} This works by an obvious modification of the preceding argument for the HLH case, and in particular Remark \ref{P:Remark1} remains valid (now the conditions~\eqref{P:207},~\eqref{P:209},~\eqref{P:210} and \eqref{P:204} are not needed, with the corresponding changes to the exceptional rules).

\subsection{The LHH case}\label{P:260}

Here we assume $N_0 \le N_1 \sim N_2$, so now \eqref{N:50} holds with
\begin{equation}\label{P:261}
  C \sim N_0^{\frac12} L_1^{\frac12}
  \left[\min\left(
  N_0^2, N_1L_2
  \right)\right]^{\frac12}.
\end{equation}
To resolve the minimum, we split into the cases $L_2 \le N_0^2/N_1$ and $L_2 > N_0^2/N_1$.

\subsubsection{The subcase $L_2 \le N_0^2/N_1$}\label{P:262:0} We restrict accordingly in \eqref{P:222}. Since also $L_2 \ge 1$, we must assume $N_0 \ge N_1^{\frac12}$. Now $C^2 \sim N_0N_1L_1L_2$, hence
\begin{equation}\label{P:262}
  S_{\Nbold}
  \lesssim
  N_0^{\frac12} N_1^{\frac12}
  \sigma_{\frac12}\left(\frac{N_0^2}{N_1}\right)
  \bignorm{F_0^{N_0}}
  \bignorm{F_1^{N_1}}
  \bignorm{F_2^{N_2}},
\end{equation}
where $\sigma_p$ is defined as in~\eqref{P:236}. Applying~\eqref{P:238} we get, for some $A,B \in \R$ depending on the $b$'s,
\begin{equation}\label{P:266}
\begin{aligned}
  \sum_{\Nbold}
  \frac{\chi_{N_1^{\frac12} \le N_0 \le N_1 \sim N_2}S_{\Nbold}}
  {N_0^{s_0}N_1^{s_1}N_2^{s_2}}
  &\lesssim
  \sum_{\Nbold}
  \chi_{N_1^{\frac12} \le N_0 \le N_1 \sim N_2}
  \frac{N_0^{A}}
  {N_1^{B}}
  \bignorm{F_0^{N_0}}
  \bignorm{F_1^{N_1}}
  \bignorm{F_2^{N_2}}
  \\
  &\lesssim
  \norm{F_0}
  \sum_{N_1,N_2} \chi_{N_1 \sim N_2}
  \frac{\Xi_A(N_1)}
  {N_1^{B}}
  \bignorm{F_1^{N_1}}
  \bignorm{F_2^{N_2}},
\end{aligned}
\end{equation}
where
$$
  \Xi_A(N)
  =
  \sum_{M} \chi_{N^{\frac12} \le M \lesssim N} M^{A}
  \sim
  \begin{cases}
  N^{A} &\text{if $A > 0$},
  \\
  \log\angles{N} &\text{if $A = 0$},
  \\
  N^{A/2} &\text{if $A < 0$}.
  \end{cases}
$$
The desired estimate follows provided that (i) $B \ge A$, (ii) $2B \ge A$ and (iii) we exclude $A=B=0$, since this guarantees $\Xi_A(N_1) \lesssim N_1^B$, hence we can apply the Cauchy-Schwarz inequality. Logarithmic factors are estimated as in \eqref{N:38}.

\begin{itemize}
\item
If the first alternative in \eqref{P:238} holds, then $A=\frac52-2(b_1+b_2)-s_0$ and $B=s_1+s_2+\frac12-(b_1+b_2)$, so $B \ge A$ and $2B \ge A$ are the same as~\eqref{P:203} and~\eqref{P:207}, respectively. Moreover, if $A=0$, then~\eqref{P:216} guarantees that $B > 0$.
\item
In the cases where we pick up the bound $1$ in \eqref{P:238}, then $A=\frac12-s_0$ and $B=s_1+s_2-\frac12$, so $B \ge A$ and $2B \ge A$ are the same as~\eqref{P:206} and~\eqref{P:207}, respectively, and $A=0$ implies $B > 0$, in view of~\eqref{P:216}.
\item
If the third to last alternative in \eqref{P:238} applies, then $A=\frac32-2b_2-s_0$ and $B=s_1+s_2-b_2$, so $B \ge A$ and $2B \ge A$ are the same as~\eqref{P:205} and~\eqref{P:207}, respectively, and $A=0$ implies $B > 0$, by~\eqref{P:216}.
\item
In the cases where we bound by one or two logarithmic factors alone, $A=\frac12+2\varepsilon-s_0$ and $B=s_1+s_2+\varepsilon-\frac12$. But this only comes up if either $b_2=\frac12$ or $b_1+b_2=1$, hence~\eqref{P:206} is strict, by~\eqref{P:212b} or~\eqref{P:215}, respectively. Therefore, $B > A$ for $\varepsilon > 0$ small enough, whereas $2B \ge A$ is the same as~\eqref{P:207}.
\item
Finally, if the fourth alternative in \eqref{P:238} prevails, then $A=\frac32-2b_2+2\varepsilon-s_0$ and $B=s_1+s_2-b_2+\varepsilon$. But now $b_1 = \frac12$, so~\eqref{P:211a} guarantees that~\eqref{P:205} is strict, hence $B > A$ for $\varepsilon > 0$ small enough, whereas $2B \ge A$ is the same as~\eqref{P:207}.
\end{itemize}

\subsubsection{The subcase $L_2 > N_0^2/N_1$}\label{P:269} Applying this restriction in \eqref{P:222}, and noting that \eqref{N:50} now holds with $C^2 \sim N_0^3L_1$, we get
\begin{equation}\label{P:270}
  S_{\Nbold}
  \lesssim
  N_0^{\frac32}
  \gamma\left(1+\frac{N_0^2}{N_1}\right)
  \bignorm{F_0^{N_0}}
  \bignorm{F_1^{N_1}}
  \bignorm{F_2^{N_2}},
\end{equation}
where $\gamma$ is defined as in~\eqref{P:244}. Since $\gamma$ is evaluated at $1+N_0^2/N_1$, we split further into the subcases $N_0 \le N_1^{\frac12}$ and $N_0 > N_1^{\frac12}$.

First, if $N_0 \le N_1^{\frac12}$, then by \eqref{P:246} the $\gamma$-factor in \eqref{P:270} is $O(1)$, hence
\begin{equation}\label{P:272}
  \sum_{\Nbold}
  \frac{\chi_{N_1 \sim N_2} \chi_{N_0 \le N_1^{\frac12}} S_{\Nbold}}
  {N_0^{s_0}N_1^{s_1}N_2^{s_2}}
  \lesssim
  \norm{F_0}
  \sum_{N_1,N_2}
  \chi_{N_1 \sim N_2}
  \frac{\Sigma_A(N_1^{\frac12})}
  {N_1^B} 
  \bignorm{F_1^{N_1}}
  \bignorm{F_2^{N_2}},
\end{equation}
where $A=\frac32 - s_0$, $B=s_1+s_2$ and $\Sigma_A$ is defined as in \eqref{D:42}. Thus, it suffices to check that (i) $2B \ge A$, (ii) $B \ge 0$ and (iii) we exclude $A=B=0$. But this follows from, respectively, \eqref{P:207}, \eqref{P:208} and \eqref{P:216}.

Second, if $N_0 > N_1^{\frac12}$, then $\gamma\left(1+\frac{N_0^2}{N_1}\right)
  \sim
  \gamma\left(\frac{N_0^2}{N_1}\right)$, so applying \eqref{P:246} we get~\eqref{P:266} for some $A,B \in \R$, and then it is enough to check that (i) $B \ge A$, (ii) $2B \ge A$ and (iii) we exclude $A=B=0$.

\begin{itemize}
\item
If $b_1 < \frac12$, then $A=\frac52-2(b_1+b_2)-s_0$ and $B=s_1+s_2+\frac12-(b_1+b_2)$, so $B \ge A$ and $2B \ge A$ are the same as~\eqref{P:203} and~\eqref{P:207}, respectively. Moreover, $A=0$ implies $B > 0$, in view of \eqref{P:216}.
\item
If $b_1 > \frac12$, then $A=\frac32-2b_2-s_0$ and $B=s_1+s_2-b_2$, so $B \ge A$ and $2B \ge A$ are the same as~\eqref{P:205} and~\eqref{P:207}, respectively, and $A=0$ implies $B > 0$, by \eqref{P:216}.
\item
If $b_1 = \frac12$, then $A=\frac32-2b_2+2\varepsilon-s_0$ and $B=s_1+s_2-b_2+\varepsilon$ for any $\varepsilon > 0$. But \eqref{P:211a} guarantees that $B > A$ for $\varepsilon > 0$ small enough, and $2B \ge A$ is the same as~\eqref{P:207}.
\end{itemize}

\subsection{The case $N_0 \ll N_1 \sim N_2$ with $\pm_1=\pm_2$}\label{P:280}

For later use we observe that the hypotheses can then be relaxed. Taking into account~\eqref{N:56}, we conclude that~\eqref{N:50} now holds with
\begin{equation}\label{P:281}
  C \sim N_0 L_1^{\frac12}
  \left[\min\left(
  N_0, L_2
  \right)\right]^{\frac12},
\end{equation}
hence the argument in \S\ref{P:230} applies (after a relabeling of the $N's$), and in particular Remark \ref{P:Remark1} applies.

\subsection{The case $b_1+b_2=\frac12$}\label{P:300} Then \eqref{P:203} becomes $s_0+s_1+s_2 \ge \frac32$, and \eqref{P:204}--\eqref{P:207} are redundant. The argument in \S\ref{P:100} does not quite work, since we only get \eqref{P:110} up to multiplication of the right hand side by the divergent sum
$$
  \sum_{L_2} L_2^{-b_2} \sum_{L_1 \le L_2} L_1^{\frac12-b_1}
  \sim \sum_{L_2} L_2^{-b_2} L_2^{\frac12-b_1}
  =
  \sum_{L_2} 1.
$$
So we must find ways to avoid this divergence. There is no problem if we restrict to $L_1 \sim L_2$, since then $
  S_{\Nbold}
  \lesssim
  \bigl(\Nmin^{012}\bigr)^{\frac32}
  \bignorm{F_0^{N_0}}
  \sum_{L_1 \sim L_2}
  \bignorm{F_1^{N_1,L_1}}
  \bignorm{F_2^{N_2,L_2}},
$
and \eqref{P:114} follows by the Cauchy-Schwarz inequality.

So from now on we restrict the summation in $S_{\Nbold}$ to $L_1 \ll L_2$. We also split $S_{\Nbold}$ depending on the signs $(\pm_1,\pm_2)$ of the temporal frequencies $(\tau_1,\tau_2)$, as in \eqref{A:50}. It is enough to estimate $S_{\Nbold}^{(+,+)}$ and $S_{\Nbold}^{(+,-)}$. We shall use the fact that, since $b_1 \in (0,\frac12)$ and $b_1+b_2=\frac12$, \eqref{P:238} gives
\begin{equation}\label{P:312}
  \sigma_p(N) \lesssim N^p,
\end{equation}
for all $p > 0$. 

\subsubsection{The case $(+,+)$}

Then \eqref{N:50} holds with
\begin{equation}\label{P:311}
  C
  \sim
  \Nmin^{012} L_1^{\frac12} 
  \left[\min\left( \Nmin^{012}, L_2 \right)\right]^{\frac12},
\end{equation}
and to resolve the minimum we split into $L_2 \lesssim \Nmin^{012}$ and $L_2 \gg \Nmin^{012}$.

If $L_2 \lesssim \Nmin^{012}$, then
$
  S_{\Nbold}^{(+,+)}
  \lesssim
  \Nmin^{012}
  \sigma_{\frac12}\bigl(\Nmin^{012}\bigr)
  \bignorm{F_0^{N_0}}
  \bignorm{F_1^{N_1}}
  \bignorm{F_2^{N_2}}
$
with $\sigma_p$ defined as in \eqref{P:236}, and \eqref{P:312} implies \eqref{P:114}, so we are done.

Now assume $L_2 \gg \Nmin^{012}$. Since $\tau_0+\tau_1+\tau_2 = 0$ in $J(\cdots)$, we have
\begin{equation}\label{P:320}
  \left(-\tau_0-\abs{\xi_0}\right)
  + \left(-\tau_1+\abs{\xi_1}\right)
  + \left(-\tau_2+\abs{\xi_2}\right)
  - \left(\abs{\xi_1}+\abs{\xi_2}-\abs{\xi_0}\right) = 0,
\end{equation}
so in absolute value, the two largest of the four terms in parentheses must be comparable. But the second term is $O(L_1)$, so it is negligible compared to the third term, which is comparable to $L_2$ in absolute value. As for the fourth term, its absolute value is comparable to, since $\xi_0+\xi_1+\xi_2 = 0$,
$$
  \min\left(\abs{\xi_1},\abs{\xi_2}\right)\theta(\xi_1,\xi_2)^2
  \sim \Nmin^{12} \theta(\xi_1,\xi_2)^2,
$$
which is negligible compared to the third term if $N_0 \sim \Nmax^{12}$, since then $L_2 \gg \Nmin^{12}$. So if $N_0 \sim \Nmax^{12}$, then $\bigabs{\tau_0+\abs{\xi_0}}
  \sim
  \bigabs{-\tau_2+\abs{\xi_2}} \sim L_2$, and \eqref{N:50} gives $
  S_{\Nbold}^{(+,+)}
  \lesssim
  \bigl(\Nmin^{12}\bigr)^{\frac32}
  \bignorm{F_1^{N_1}}
  \sum_{L_2}
  \bignorm{\chi_{\abs{\tau_0+\abs{\xi_0}} \sim L_2} F_0^{N_0}}
  \bignorm{F_2^{N_2,L_2}},
$
so we can sum $L_2$ using the Cauchy-Schwarz inequality, thus obtaining \eqref{P:114}.

It remains to consider the case
$$
  N_0 \ll N_1 \sim N_2.
$$
Then the preceding argument does not work, since we only know that $L_2 \gg N_0$, whereas the fourth term in \eqref{P:320} is comparable to $N_1 \gg N_0$. To get around this problem, we apply first a standard decomposition into cubes. Letting $\mathcal Q$ be a tiling of $\R^3$ into almost disjoint cubes $Q$ of sidelength $N_0$, we write
$$
  S_{\Nbold}^{(+,+)}
  =
  \sum_{Q_1,Q_2 \in \mathcal Q} a_{\Nbold}^{(Q_1,Q_2)},
$$
where
\begin{multline*}
  a_{\Nbold}^{(Q_1,Q_2)}
  =
  \sum_{L_1,L_2}
  \chi_{L_1 \ll L_2}
  \chi_{N_0 \ll L_2} \chi_{N_0 \ll N_1 \sim N_2}
  L_1^{-b_1}L_2^{-b_2}
  \\
  \times
  \triple{F_0^{N_0}}{\chi_{\R \times Q_1} F_1^{N_1,L_1,+}}{\chi_{\R \times Q_2} F_2^{N_2,L_2,+}}.
\end{multline*}
In the integral defining $J(\cdots)$, $\xi_0$ is now restricted to the ball $\setb{ \xi }{ \abs{\xi} \le cN_0 }$ for some absolute constant $c \ge 1$. On the other hand, $\xi_0=-\xi_1-\xi_2$ where $\xi_1 \in Q_1$ and $\xi_2 \in Q_2$. Therefore, once $Q_1 \in \mathcal Q$ has been chosen, the choice of $Q_2$ is limited to a subset $\mathcal Q(Q_1) \subset \mathcal Q$ of cardinality $O(1)$. Thus, it will be enough to show
\begin{equation}\label{P:340}
  a_{\Nbold}^{(Q_1,Q_2)}
  \lesssim
  N_0^{\frac32}
  \bignorm{F_0^{N_0}}
  \norm{\chi_{\R \times Q_1}F_1^{N_1}}
  \norm{\chi_{\R \times Q_2}F_2^{N_2}},
\end{equation}
since then we can just sum over $Q_1 \in \mathcal Q$ and $Q_2 \in \mathcal Q(Q_1)$ and apply the Cauchy-Schwarz inequality to obtain the corresponding inequality for $S_{\Nbold}^{(+,+)}$.

We have
\begin{equation}\label{P:344}
  Q_1 = \xi_1^* + [0,N_0]^3, \qquad
  Q_2 = \xi_2^* + [0,N_0]^3,
\end{equation}
for some $\xi_1^*,\xi_2^* \in \R^3$ such that
\begin{equation}\label{P:346}
  \abs{\xi_1^*}, \abs{\xi_2^*} \sim N_1 \sim N_2 \gg N_0,
  \qquad
  \abs{\xi_1^* + \xi_2^*} \lesssim N_0.
\end{equation}
From~\eqref{P:320} we get
\begin{equation}\label{P:347}
\begin{aligned}
  \tau_0 + \abs{\xi_0}
  &=
  (- \tau_1 + \abs{\xi_1})
  +
  (- \tau_2 + \abs{\xi_2})
  +
  \abs{\xi_1}+\abs{\xi_2}-\abs{\xi_1+\xi_2}
  \\
  &=
  O(L_1) + (- \tau_2 + \abs{\xi_2}) + \abs{\xi_1^*}+\abs{\xi_2^*}-\abs{\xi_1^* + \xi_2^*} + O(N_0).
\end{aligned}
\end{equation}
But $\bigabs{- \tau_2 + \abs{\xi_2}} \sim L_2$, whereas $L_1 \ll L_2$ and $N_0 \ll L_2$. We conclude that
$$
  \tau_0 + \abs{\xi_0}
  \in p + I_{L_2},
$$
where $p
  =
  p(Q_1,Q_2)
  =
  \abs{\xi_1^*}+\abs{\xi_2^*}-\abs{\xi_1^* + \xi_2^*}$ and
  $I_{L_2}
  =
  \left[-cL_2,-c^{-1}L_2\right] \cup \left[c^{-1}L_2,cL_2\right]$
for some absolute constant $c \gg 1$.

By the Sobolev type estimate,~\eqref{N:42}, we can therefore dominate $a_{\Nbold}^{(Q_1,Q_2)}$ by
$$
  N_0^{\frac32}
  \norm{\chi_{\R \times Q_1} F_1^{N_1}}
  \left(\sum_{L_2}
  \norm{\chi_{\tau_0 + \abs{\xi_0} \in p + I_{L_2}} F_0^{N_0}}
  \norm{\chi_{\R \times Q_2} F_2^{N_2,L_2,+}}\right),
$$
and~\eqref{P:340} then follows by the Cauchy-Schwarz inequality, since
$$
  \left(\sum_{L_2}
  \norm{\chi_{\tau_0 + \abs{\xi_0} \in p + I_{L_2}} F_0^{N_0}}^2 \right)^{\frac12}
  =
  \norm{\left(\sum_{L_2}\chi_{\tau_0 + \abs{\xi_0} \in p + I_{L_2}}\right) F_0^{N_0}}
  \sim
  \bignorm{F_0^{N_0}}.
$$

\subsubsection{The case $(+,-)$}\label{P:348}

First observe that when $N_0 \sim \Nmax^{12}$ we have exactly the same estimates as in the $(+,+)$ case, hence we proceed as we did there. The only difference is that \eqref{P:320} must now be replaced by one of the following two:
\begin{align}
  \label{P:350}
  \left(-\tau_0+\abs{\xi_0}\right)
  + \left(-\tau_1+\abs{\xi_1}\right)
  + \left(-\tau_2-\abs{\xi_2}\right)
  - \left(\abs{\xi_0} + \abs{\xi_1} - \abs{\xi_2}\right) &= 0,
  \\
  \label{P:352}
  \left(-\tau_0-\abs{\xi_0}\right)
  + \left(-\tau_1+\abs{\xi_1}\right)
  + \left(-\tau_2-\abs{\xi_2}\right)
  + \left(\abs{\xi_0} + \abs{\xi_2} - \abs{\xi_1}\right) &= 0.
\end{align}
Specifically, we use \eqref{P:350} if $\abs{\xi_1} \le \abs{\xi_2}$ and \eqref{P:352} otherwise. Then the fourth terms in~\eqref{P:350} or \eqref{P:352}, respectively, are dominated in absolute value by $\Nmin^{012}$, whereas in the $(+,+)$ case we had the bound $\Nmin^{12}$.

We are then left with the case
$$
  N_0 \ll N_1 \sim N_2.
$$
The estimates are then less favorable than in the $(+,+)$ case, since now \eqref{N:50} holds with $C$ as in \eqref{P:261}.

If we restrict to $L_2 \lesssim N_0^2/N_1$, then
$$
  S_{\Nbold}^{(+,-)}
  \lesssim
  N_0^{\frac12} N_1^{\frac12}
  \sigma_{\frac12}\left( \frac{N_0^2}{N_1} \right)
  \bignorm{F_0^{N_0}}
  \bignorm{F_1^{N_1}}
  \bignorm{F_2^{N_2}},
$$
and in view of \eqref{P:312} we then get the desired estimate.

Now consider $L_2 \gg N_0^2/N_1$. We reduce to proving \eqref{P:340} for
\begin{multline*}
  a_{\Nbold}^{(Q_1,Q_2)}
  =
  \sum_{L_1,L_2}
  \chi_{L_1 \ll L_2}
  \chi_{\frac{N_0^2}{N_1} \ll L_2}
  \chi_{N_0 \ll N_1 \sim N_2}
  L_1^{-b_1}L_2^{-b_2}
  \\
  \times
  \triple{F_0^{N_0}}{\chi_{\R \times Q_1} F_1^{N_1,L_1,+}}{\chi_{\R \times Q_2} F_2^{N_2,L_2,-}},
\end{multline*}
where $Q_1$ and $Q_2$ are as in \eqref{P:344}--\eqref{P:346}.

If we use \eqref{N:42}, we get \eqref{P:340} up to multiplication by the sum (recall $b_1 \in (0,\frac12)$ and $b_1+b_2=\frac12$)
$$
  \sum_{L_1,L_2}
  \chi_{L_1 \ll L_2}
  \chi_{\frac{N_0^2}{N_1} \ll L_2}
  L_1^{\frac12-b_1} L_2^{-b_2}
  \sim
  \sum_{L_2}
  \chi_{\frac{N_0^2}{N_1} \ll L_2}
  1
$$
but this diverges, of course. To avoid this divergence, we shall use some orthogonality properties.

Since we are in the case $N_0 \ll N_1 \sim N_2$ with opposite signs, it makes sense to decompose using thickened null hyperplanes instead of cones. So as a replacement for~\eqref{P:350}--\eqref{P:352} we try the following: Set $\omega = \xi_1^*/\abs{\xi_1^*}$, where $\xi_1^*$ is one corner of the cube $Q_1$, as in \eqref{P:344}--\eqref{P:346}. Since $\tau_0+\tau_1+\tau_2 = 0$ and $\xi_0+\xi_1+\xi_2 = 0$, we have
\begin{equation}\label{P:360}
  (-\tau_0 + \xi_0 \cdot \omega)
  +
  (-\tau_1 + \abs{\xi_1})
  +
  (-\tau_2 - \abs{\xi_2})
  -
  (\abs{\xi_1} - \xi_1 \cdot \omega)
  +
  (\abs{\xi_2} + \xi_2 \cdot \omega)
  =
  0.
\end{equation}
Note that
\begin{equation}\label{P:362}
  \abs{\xi_1} - \xi_1 \cdot \omega
  =
  O\left(N_1 \theta(\xi_1,\omega)^2\right)
  =
  O\left(\frac{N_0^2}{N_1}\right).
\end{equation}
where the last equality holds since $\xi_1 \in Q_1$, hence $\theta(\xi_1,\omega) \lesssim N_0/N_1$.

Similarly, since $\xi_2 \in Q_2$, and since $-Q_2$ is within an $O(N_0)$-neighborhood of $Q_1$, by the assumption \eqref{P:344}, we find that
\begin{equation}\label{P:364}
  \abs{\xi_2} + \xi_2 \cdot \omega
  =
  O\left(\frac{N_0^2}{N_1}\right).
\end{equation}
Since $N_0^2/N_1 \ll L_2$ and $L_1 \ll L_2$, we can conclude from \eqref{P:360}--\eqref{P:364} that
$$
  \abs{-\tau_0 + \xi_0 \cdot \omega}
  \sim
  \bigabs{-\tau_2 - \abs{\xi_2}}
  \sim L_2,
$$
hence $L_2$ can be summed by the Cauchy-Schwarz inequality.

\section{The case $0 < b_0,b_1,b_2$}\label{P:398}

Then the product law reads:

\begin{theorem}\label{P:Thm3} Let $n = 3$. Assume
\begin{align}
  \label{P:400}
  &b_0, b_1, b_2 > 0
  \\
  \label{P:401}
  &b_0 + b_1 + b_2 \ge \frac12
  \\
  \label{P:402}
  &s_0 + s_1 + s_2 \ge 2 - (b_0 + b_1 + b_2)
  \\
  \label{P:403}
  &s_0 + s_1 + s_2 \ge \frac32 - (b_0 + b_1)
  \\
  \label{P:404}
  &s_0 + s_1 + s_2 \ge \frac32 - (b_0 + b_2)
  \\
  \label{P:405}
  &s_0 + s_1 + s_2 \ge \frac32 - (b_1 + b_2)
  \\
  \label{P:406}
  &s_0 + s_1 + s_2 \ge 1
  \\
  \label{P:407}
  &(s_0 + b_0) + 2s_1 + 2s_2 \ge \frac32& &\Hproduct{\text{L} & \text{H} & \text{H} \\ & + & -}
  \\
  \label{P:408}
  &2s_0 + (s_1 + b_1) + 2s_2 \ge \frac32& &\Hproduct{\text{H} & \text{L} & \text{H} \\ + & & -}
  \\
  \label{P:409}
  &2s_0 + 2s_1 + (s_2 + b_2) \ge \frac32& &\Hproduct{\text{H} & \text{H} & \text{L} \\ + & - & }
  \\
  \label{P:410}
  &s_1 + s_2 \ge 0& &\Hproduct{\text{L} & \text{H} & \text{H} \\ & \phantom{+} & \phantom{-}}
  \\
  \label{P:411}
  &s_0 + s_2 \ge 0& &\Hproduct{\text{H} & \text{L} & \text{H} \\ \phantom{+} & & \phantom{-}}
  \\
  \label{P:412}
  &s_0 + s_1 \ge 0& &\Hproduct{\text{H} & \text{H} & \text{L} \\ \phantom{+} & \phantom{-} & },
\end{align}
as well as the exceptions:
\begin{align}
  \label{P:413}
  &\parbox[t]{10.8cm}{If $b_0=\frac12$, then~\eqref{P:402}=\eqref{P:405} must be strict.}
  \\
  \label{P:414}
  &\parbox[t]{10.8cm}{If $b_1=\frac12$, then~\eqref{P:402}=\eqref{P:404} must be strict.}
  \\
  \label{P:415}
  &\parbox[t]{10.8cm}{If $b_2=\frac12$, then~\eqref{P:402}=\eqref{P:403} must be strict.}
  \\
  \label{P:416}
  &\parbox[t]{10.8cm}{If $b_0+b_1=\frac12$, then~\eqref{P:403}=\eqref{P:406} must be strict.}
  \\
  \label{P:417}
  &\parbox[t]{10.8cm}{If $b_0+b_2=\frac12$, then~\eqref{P:404}=\eqref{P:406} must be strict.}
  \\
  \label{P:418}
  &\parbox[t]{10.8cm}{If $b_1+b_2=\frac12$, then~\eqref{P:405}=\eqref{P:406} must be strict.}
  \\
  \label{P:419}
  &\parbox[t]{10.8cm}{If $b_0+b_1+b_2=1$, then~\eqref{P:402}=\eqref{P:406} must be strict.}
  \\
  \label{P:420}
  &\parbox[t]{10.8cm}{We require~\eqref{P:407} to be strict if $s_0+b_0$ takes one of the values $\frac32$, $\frac12 + 2b_0$, $\frac32 - 2b_1$, $\frac32 - 2b_2$ or $\frac52 - 2(b_1+b_2)$.}
  \\
  \label{P:421}
  &\parbox[t]{10.8cm}{We require~\eqref{P:408} to be strict if $s_1+b_1$ takes one of the values $\frac32$, $\frac12 + 2b_1$, $\frac32 - 2b_0$, $\frac32 - 2b_2$ or $\frac52 - 2(b_0+b_2)$.}
  \\
  \label{P:422}
  &\parbox[t]{10.8cm}{We require~\eqref{P:409} to be strict if $s_2+b_2$ takes one of the values $\frac32$, $\frac12 + 2b_2$, $\frac32 - 2b_0$, $\frac32 - 2b_1$ or $\frac52 - 2(b_0+b_1)$.}
  \\ 
  \label{P:423}
  &\parbox[t]{10.8cm}{If one of~\eqref{P:402}--\eqref{P:406} is an equality, then~\eqref{P:410}--\eqref{P:412} must be strict.}
\end{align}
Then $\Hproduct{s_0 & s_1 & s_2 \\ b_0 & b_1 & b_2}$ is a product.
\end{theorem}

\begin{remark}\label{P:424} For later use we note that
\begin{equation}\label{P:474}
  s_0 + 2(s_1+s_2) \ge 1 + 2\varepsilon,
\end{equation}
for some $\varepsilon > 0$. This follows from~\eqref{P:406} if $s_1+s_2>0$. If $s_1+s_2=0$, on the other hand, then we infer from~\eqref{P:423} that~\eqref{P:406} is strict, so again the desired inequality holds. Applying the same argument to \eqref{P:402}, we find that
\begin{equation}\label{P:476}
  s_0 + 2(s_1+s_2) + b_0 + b_1 +b_2 \ge 2 + 2\varepsilon,
\end{equation}
for some $\varepsilon > 0$.
\end{remark}

We now prove Theorem \ref{P:Thm3}. By symmetry we may assume $L_0 \ge L_2 \ge L_1$. Then
$
  P = \Hproduct{s_0 & s_1 & s_2 \\ b_0 & b_1 & b_2}
$
is a product if
$
  P' = \Hproduct{s_0 & s_1 & s_2 \\ 0 & b_1 & b_0+b_2}
$
is. So we go ahead and check whether $P'$ satisfies the hypotheses of Theorem \ref{P:Thm2} (replace $b_2$ by $b_2'=b_0+b_2$ there). This is indeed seen to be the case if we restrict to the HLH and HHL interactions, since then Remark~\ref{P:Remark1} applies. Furthermore, the LHH interaction is also admissible if $b_0+b_1+b_2=\frac12$, since then~\eqref{P:403}--\eqref{P:409} are all strict, and in particular~\eqref{P:204}--\eqref{P:207} are strict for $P'$, hence the rules~\eqref{P:211a}--\eqref{P:216} are redundant.

In view of these reductions, we may assume
$
  b_0 + b_1 + b_2 > \frac12,
$
and we need only consider the case
$
  N_0 \ll N_1 \sim N_2.
$
Then we shall prove \eqref{P:110} with
\begin{equation}\label{P:440}
  S_{\Nbold}
  =
  \sum_{\Lbold} \chi_{L_1 \le L_2 \le L_0}
  \frac{
  \triple{F_0^{N_0,L_0}}{F_1^{N_1,L_1}}{F_2^{N_2,L_2}}
  }
  {L_0^{b_0} L_1^{b_1} L_2^{b_2}}
\end{equation}
where now $\Lbold=(L_0,L_1,L_2)$. By the estimates in \S\ref{N:40} we deduce that
\begin{equation}\label{P:460}
  \triple{F_0^{N_0,L_0}}{F_1^{N_1,L_1}}{F_2^{N_2,L_2}}
  \le C
  \bignorm{F_0^{N_0}}
  \bignorm{F_1^{N_1}}
  \bignorm{F_2^{N_2}}
\end{equation}
holds with
\begin{equation}\label{P:462}
  C
  \sim
  N_0^{\frac12} L_1^{\frac12}
  \left[ \min\left(
  N_0^2,
  N_1L_2,
  N_0L_0
  \right) \right]^{\frac12}.
\end{equation}
To resolve the minimum, we distinguish $L_0 \ge N_0$ and $L_0 < N_0$, and in the latter case we split further into $L_2 \le (N_0/N_1)L_0$ and $L_2 > (N_0/N_1)L_0$.

\subsection{The case $L_0 \ge N_0$}\label{P:464} Then we remove $N_0L_0$ from the minimum in~\eqref{P:462}, and we sum out $L_0$ using
$$
  \sum_{L_0 \ge N_0} L_0^{-b_0} \sim N_0^{-b_0},
$$
which holds since $b_0 > 0$. We can then proceed as in \S\ref{P:260}, but replacing $s_0$ there by $s_0'=s_0+b_0$. Thus, we replace
$
  P = \left.\Hproduct{s_0 & s_1 & s_2 \\ b_0 & b_1 & b_2}\right\vert_{\text{LHH}}
$
by
$
  P' = \left.\Hproduct{s_0+b_0 & s_1 & s_2 \\ 0 & b_1 & b_2}\right\vert_{\text{LHH}}.
$
Again we go ahead and check whether the hypotheses on $P$ imply the relevant conditions on $P'$ in \S\ref{P:200}. This is indeed seen to be the case if $b_1+b_2 \ge \frac12$. The only point which is not completely trivial is that the rule \eqref{P:420} takes care of all the exceptional values in \eqref{P:216} apart from $s_0'=s_0+b_0=\frac12$. But if $s_0+b_0=\frac12$, then~\eqref{P:407} (which is the same as~\eqref{P:207} for $P'$) must be strict, for if it were an equality we would have $s_1+s_2 = \frac12$, but then \eqref{P:406} implies $s_0 \ge \frac12$, contradicting $s_0+b_0=\frac12$. Thus, $P'$ is indeed a product if $b_1 + b_2 \ge \frac12$.

This still leaves the case
$$
  b_1+b_2 < \frac12.
$$
Now we do not sum $L_0$ right away, but repeat instead the LHH argument in \S\ref{P:260} as far as possible; the argument only fails because we use~\eqref{P:246} to estimate the $\gamma$-factor in \eqref{P:270}, but now the sum in~\eqref{P:246} diverges, since $b_1+b_2 < \frac12$. However, we now have $L_2 \le L_0$, so the divergent sum can be replaced by
$$
  \delta(L_0)
  =
  \sum_{L_1,L_2}
  \chi_{L_1 \le L_2 \le L_0} L_1^{\frac12-b_1} L_2^{-b_2}
  \sim
  \sum_{L_2 \le L_0} L_2^{\frac12-b_1-b_2}
  \sim L_0^{\frac12-b_1-b_2}.
$$
Thus, $\gamma(\ldots)$ in \eqref{P:270} can be replaced by
$
  \sum_{L_0 \ge N_0} L_0^{\frac12-b_0-b_1-b_2} \sim N_0^{\frac12-b_0-b_1-b_2},
$
hence
$$
  S_{\Nbold}
  \lesssim
  N_0^{2-b_0-b_1-b_2}
  \bignorm{F_0^{N_0}}
  \bignorm{F_1^{N_1}}
  \bignorm{F_2^{N_2}},
$$
so setting $A=2-s_0-b_0-b_1-b_2$ and $B=s_1+s_2$, it suffices to check that (i) $B \ge A$, (ii) $B \ge 0$ and (iii) we exclude $A=B=0$. But this follows from, respectively, \eqref{P:402}, \eqref{P:410} and \eqref{P:423}.

\subsection{The case $L_0 < N_0$ with $L_2 \le (N_0/N_1)L_0$}\label{P:466}  Then $N_1 \le N_1L_2 \le N_0L_0 \le N_0^2$, hence $N_1^{\frac12} \le N_0$. Since \eqref{P:462} now reads $C^2 \sim N_0N_1L_1L_2$, we get
$$
  S_{\Nbold}
  \lesssim
  N_0^{\frac12}N_1^{\frac12}
  \kappa_{\frac12}\left(N_0,\frac{N_0}{N_1}\right)
  \bignorm{F_0^{N_0}}
  \bignorm{F_1^{N_1}}
  \bignorm{F_2^{N_2}},
$$
where we write, for $0 < r \le 1$ and $p > 0$,
\begin{align*}
  \kappa_p(N_0,r)
  &=
  \sum_{\Lbold} \chi_{L_1 \le L_2 \le r L_0}
  \chi_{L_0 \le N_0}
  L_0^{-b_0} L_1^{\frac12-b_1} L_2^{p-b_2}
  \\
  &=
  \sum_{L_0} \chi_{r^{-1} \le L_0 \le N_0}
  L_0^{-b_0} \sigma_p(rL_0),
\end{align*}
with $\sigma_p$ as in~\eqref{P:236}. Then by~\eqref{P:238},
$$
  \kappa_p(N_0,r)
  \lesssim
  r^{b_0}
  \sum_{L_0} \chi_{1 \le rL_0 \le rN_0}
  \begin{cases}
  (rL_0)^{\frac12+p-b_0-b_1-b_2} &\text{$b_1 < \frac12$, $b_1+b_2 < \frac12+p$}
  \\
  (rL_0)^{-b_0} \log\angles{rL_0} &\text{$b_1 < \frac12$, $b_1+b_2 = \frac12+p$}
  \\
  (rL_0)^{-b_0}  &\text{$b_1 < \frac12$, $b_1+b_2 > \frac12+p$}
  \\
  (rL_0)^{p-b_0-b_2}\log\angles{rL_0} &\text{$b_1 = \frac12$, $b_2 < p$}
  \\
  (rL_0)^{-b_0} \log^2\angles{rL_0} &\text{$b_1 = \frac12$, $b_2 = p$}
  \\
  (rL_0)^{-b_0}  &\text{$b_1 = \frac12$, $b_2 > p$}
  \\
  (rL_0)^{p-b_0-b_2} &\text{$b_1 > \frac12$, $b_2 < p$}
  \\
  (rL_0)^{-b_0} \log\angles{rL_0} &\text{$b_1 > \frac12$, $b_2 = p$}
  \\
  (rL_0)^{-b_0}  &\text{$b_1 > \frac12$, $b_2 > p$}.
  \end{cases}
$$
Thus, recalling also that $b_0 > 0$, we find that $\kappa_p=\kappa_p(N_0,r)$ satisfies
$$
  \kappa_p
  \lesssim
  \begin{cases}
  r^{b_0}(rN_0)^{\frac12+p-b_0-b_1-b_2} &\text{$b_1 < \frac12$, $b_1+b_2 < \frac12+p$, $b_0+b_1+b_2 < \frac12+p$}
  \\
  r^{b_0}\log\angles{rN_0} &\text{$b_1 < \frac12$, $b_1+b_2 < \frac12+p$, $b_0+b_1+b_2 = \frac12+p$}
  \\
  r^{b_0} &\text{$b_1 < \frac12$, $b_1+b_2 < \frac12+p$, $b_0+b_1+b_2 > \frac12+p$}
  \\
  r^{b_0} &\text{$b_1 < \frac12$, $b_1+b_2 \ge \frac12+p$}
  \\
  r^{b_0}(rN_0)^{p-b_0-b_2}\log\angles{rN_0} &\text{$b_1 = \frac12$, $b_2 < p$, $b_0+b_2 < p$}
  \\
  r^{b_0}\log^2\angles{rN_0} &\text{$b_1 = \frac12$, $b_2 < p$, $b_0+b_2 = p$}
  \\
  r^{b_0} &\text{$b_1 = \frac12$, $b_2 < p$, $b_0+b_2 > p$}
  \\
  r^{b_0} &\text{$b_1 = \frac12$, $b_2 \ge p$}
  \\
  r^{b_0}(rN_0)^{p-b_0-b_2} &\text{$b_1 > \frac12$, $b_2 < p$, $b_0+b_2 < p$}
  \\
  r^{b_0}\log\angles{rN_0} &\text{$b_1 > \frac12$, $b_2 < p$, $b_0+b_2 = p$}
  \\
  r^{b_0} &\text{$b_1 > \frac12$, $b_2 < p$, $b_0+b_2 > p$}
  \\
  r^{b_0} &\text{$b_1 > \frac12$, $b_2 \ge p$}.
  \end{cases}
$$
Applying this with $r=N_0/N_1$ and $p=\frac12$, we get \eqref{P:266} for some $A,B \in \R$, so it is enough to check that (i) $B \ge A$, (ii) $2B \ge A$ and (iii) we exclude $A=B=0$. Logarithmic factors are estimated as in \eqref{N:38}.

\begin{itemize}
\item
If $\kappa_{\frac12} \lesssim r^{b_0} = \frac{N_0^{b_0}}{N_1^{b_0}}$, then $A = \frac12-s_0+b_0$ and $B=s_1+s_2-\frac12+b_0$, so $B \ge A$ and $2B \ge A$ are the same as~\eqref{P:406} and~\eqref{P:407}, respectively. Moreover, $A=0$ implies $B > 0$, in view of \eqref{P:420}.
\item
If $\kappa_{\frac12} \lesssim r^{b_0}(rN_0)^{1-b_0-b_1-b_2}$, then $A = \frac52-s_0-b_0-2b_1-2b_2$ and $B=s_1+s_2+\frac12-b_1-b_2$, so $B \ge A$ and $2B \ge A$ are the same as~\eqref{P:402} and~\eqref{P:407}, respectively, and $A=0$ implies $B > 0$, by \eqref{P:420}.
\item
If $\kappa_{\frac12} \lesssim r^{b_0}(rN_0)^{\frac12-b_0-b_2}$, then $A = \frac32-s_0-b_0-2b_2$ and $B=s_1+s_2-b_2$, so $B \ge A$ and $2B \ge A$ are the same as~\eqref{P:404} and~\eqref{P:407}, respectively, and $A=0$ implies $B > 0$, by \eqref{P:420}.
\item
If $\kappa_{\frac12} \lesssim r^{b_0} (rN_0)^\varepsilon$, then $A=\frac12-s_0+b_0+2\varepsilon$ and $B=s_1+s_2-\frac12+b_0+\varepsilon$. But now either $b_0+b_2=\frac12$ or $b_0+b_1+b_2=1$, so by~\eqref{P:417} or~\eqref{P:419}, respectively, we have $B > A$. Moreover, $2B \ge A$ is the same as~\eqref{P:407}.
\item
Finally, if $\kappa_{\frac12} \lesssim r^{b_0} (rN_0)^{\frac12-b_0-b_2+\varepsilon}$, then $A=\frac32-s_0-b_0-2b_2+2\varepsilon$ and $B=s_1+s_2-b_2+\varepsilon$. But now $b_1=\frac12$, so rule~\eqref{P:414} implies $B > A$, whereas $2B \ge A$ again is the same as~\eqref{P:407}.
\end{itemize}

\subsection{The case $L_0 < N_0$ with $L_2 > (N_0/N_1)L_0$}

Then we see that \eqref{P:462} simplifies to $C^2 \sim N_0^2 L_0 L_1$, so
\begin{equation}\label{P:470}
  S_{\Nbold}
  \lesssim
  N_0
  \rho_{\frac12}\left(N_0,\frac{N_0}{N_1}\right)
  \bignorm{F_0^{N_0}}
  \bignorm{F_1^{N_1}}
  \bignorm{F_2^{N_2}},
\end{equation}
where we use the notation, for $0 < r \le 1$ and $p > 0$,
\begin{equation}\label{P:471}
\begin{aligned}
  \rho_p(N_0,r)
  &=
  \sum_{\Lbold} \chi_{L_1 \le L_2} \chi_{rL_0 \le L_2 \le L_0 \le N_0}
  L_0^{p-b_0} L_1^{\frac12-b_1} L_2^{-b_2}
  \\
  &=
  \sum_{L_0} \chi_{L_0 \le N_0} L_0^{p-b_0} \Gamma\bigl(\max(1,rL_0),L_0\bigr),
\end{aligned}
\end{equation}
and we write, for $1 \le A < B$,
$$
  \Gamma(A,B)
  =
  \sum_{L_1,L_2}
  \chi_{L_1 \le L_2} \chi_{A \le L_2 \le B}
  L_1^{\frac12-b_1} L_2^{-b_2}.
$$
Recalling that $b_2 > 0$, we find
\begin{equation}\label{P:472}
\begin{aligned}
  \Gamma(A,B)
  &\sim
  \sum_{A \le L_2 \le B}
  L_2^{-b_2}
  \times
  \begin{cases}
  L_2^{\frac12-b_1} &\text{if $b_1 < \frac12$}
  \\
  \log\angles{L_2} &\text{if $b_1 = \frac12$}
  \\
  1 &\text{if $b_1 > \frac12$}
  \end{cases}
  \\
  &\lesssim
  \begin{cases}
  B^{\frac12-b_1-b_2} &\text{if $b_1 < \frac12$, $b_1+b_2 < \frac12$}
  \\
  \left(\frac{B}{A}\right)^\varepsilon &\text{if $b_1 < \frac12$, $b_1+b_2 = \frac12$}
  \\
  A^{\frac12-b_1-b_2} &\text{if $b_1 < \frac12$, $b_1+b_2 > \frac12$}
  \\
  A^{\varepsilon-b_2} &\text{if $b_1 = \frac12$}
  \\
  A^{-b_2} &\text{if $b_1 > \frac12$},
  \end{cases}
\end{aligned}
\end{equation}
for any $\varepsilon > 0$.

Now we split further into $N_0 < N_1^{\frac12}$ and $N_0 \ge N_1^{\frac12}$.

\subsubsection{The subcase $N_0 < N_1^{\frac12}$} Setting $r=N_0/N_1$, we then have $rN_0 < 1$, hence
$$
  \rho_p(N_0,r) = \sum_{L_0} \chi_{L_0 \le N_0} L_0^{p-b_0} \Gamma(1,L_0),
$$
so by~\eqref{P:472} we get
\begin{align*}
  \rho_p(N_0,r)
  &\lesssim
  \sum_{L_0 \le N_0} L_0^{p-b_0}
  \times
  \begin{cases}
  L_0^{\frac12-b_1-b_2} &\text{if $b_1 < \frac12$, $b_1+b_2 < \frac12$}
  \\
  L_0^\varepsilon &\text{if $b_1 < \frac12$, $b_1+b_2 = \frac12$}
  \\
  1 &\text{if $b_1 + b_2 > \frac12$}
  \end{cases}
  \\
  &\lesssim
  \begin{cases}
  N_0^{\frac12+p-b_0-b_1-b_2} &\text{if $b_1 < \frac12$, $b_1+b_2 < \frac12$, $b_0+b_1+b_2 < \frac12+p$}
  \\
  N_0^\varepsilon &\text{if $b_1 < \frac12$, $b_1+b_2 < \frac12$, $b_0+b_1+b_2 = \frac12+p$}
  \\
  1 &\text{if $b_1 < \frac12$, $b_1+b_2 < \frac12$, $b_0+b_1+b_2 > \frac12+p$}
  \\
  N_0^{p-b_0+\varepsilon} &\text{if $b_1 < \frac12$, $b_1+b_2 = \frac12$, $b_0 \le p$}
  \\
  1 &\text{if $b_1 < \frac12$, $b_1+b_2 = \frac12$, $b_0 > p$}
  \\
  N_0^{p-b_0} &\text{if $b_1 + b_2 > \frac12$, $b_0 < p$}
  \\
  N_0^\varepsilon &\text{if $b_1 + b_2 > \frac12$, $b_0 = p$}
  \\
  1 &\text{if $b_1 + b_2 > \frac12$, $b_0 > p$},
  \end{cases}
\end{align*}
for any $\varepsilon > 0$. Plugging this into~\eqref{P:470}, with $r=N_0/N_1$ and $p=\frac12$, we get~\eqref{P:272}, for some $A,B \in \R$, so it suffices to check that (i) $2B \ge A$, (ii) $B \ge 0$ and (iii) we exclude $A=B=0$. In fact, $B=s_1+s_2$ in all cases, so $B \ge 0$ is the same as~\eqref{P:410}.

\begin{itemize}
\item
If $\rho_{\frac12} \lesssim  N_0^\varepsilon$, then $A=1-s_0+\varepsilon$, and \eqref{P:474} implies $2B > A$. This also covers the cases where $\rho_{\frac12} \lesssim 1$, of course.
\item
If $\rho_{\frac12} \lesssim N_0^{1-b_0-b_1-b_2}$, then $A=2-s_0-b_0-b_1-b_2$, and \eqref{P:476} implies $2B > A$. 
\item
If $\rho_{\frac12} \lesssim N_0^{\frac12-b_0}$, then $A=\frac32-s_0-b_0$, and $2B \ge A$ is the same as~\eqref{P:407}. Moreover, $A=0$ implies $B > 0$, in view of~\eqref{P:420}.
\item
If $\rho_{\frac12} \lesssim N_0^{\frac12-b_0+\varepsilon}$, then $A=\frac32-s_0-b_0+\varepsilon$, so we want strict inequality in~\eqref{P:407}, since this implies $2B > A$ for $\varepsilon > 0$ small enough. Clearly, ~\eqref{P:407} is strict if $s_0 + b_0 > \frac32$, and in fact also if $s_0+b_0=\frac32$, in view of~\eqref{P:420}. So it remains to consider the case $s_0+b_0 < \frac32$, but then~\eqref{P:407} implies $s_1+s_2 > 0$, so adding $s_1+s_2$ to~\eqref{P:402} and using the fact that we are in the case $b_1+b_2=\frac12$, we see that~\eqref{P:407} is again strict.
\end{itemize}

\subsubsection{The subcase $N_0 \ge N_1^{\frac12}$} Then $rN_0 \ge 1$, where $r=N_0/N_1$, so by~\eqref{P:471} and~\eqref{P:472},
\begin{align*}
  \rho_p(N_0,r) &= \sum_{L_0} \chi_{L_0 < r^{-1}} L_0^{p-b_0} \Gamma(1,L_0)
  +
  \sum_{L_0} \chi_{r^{-1} \le L_0 \le N_0} L_0^{p-b_0} \Gamma(rL_0,L_0)
  \\
  &\lesssim
  \sum_{L_0 < r^{-1}} L_0^{p-b_0}
  \times
  \begin{cases}
  L_0^{\frac12-b_1-b_2} &\text{$b_1 < \frac12$, $b_1+b_2 < \frac12$}
  \\
  L_0^\varepsilon &\text{$b_1 < \frac12$, $b_1+b_2 = \frac12$}
  \\
  1 &\text{$b_1 + b_2 > \frac12$}
  \end{cases}
  \\
  &\quad +
  \sum_{r^{-1} \le L_0 \le N_0} L_0^{p-b_0}
  \times
  \begin{cases}
  L_0^{\frac12-b_1-b_2} &\text{$b_1 < \frac12$, $b_1+b_2 < \frac12$}
  \\
  r^{-\varepsilon} &\text{$b_1 < \frac12$, $b_1+b_2 = \frac12$}
  \\
  (rL_0)^{\frac12-b_1-b_2} &\text{$b_1 < \frac12$, $b_1+b_2 > \frac12$}
  \\
  (rL_0)^{\varepsilon-b_2} &\text{$b_1 = \frac12$}
  \\
  (rL_0)^{-b_2} &\text{$b_1 > \frac12$},
  \end{cases}
\end{align*}
for any $\varepsilon > 0$. From this we conclude that $\rho_p = \rho_p(N_0,r)$ verifies the estimates
$$
  \rho_p \lesssim
  \begin{cases}
  N_0^{\frac12+p-b_0-b_1-b_2} &\text{$b_1 < \frac12$, $b_1+b_2 < \frac12$, $b_0+b_1+b_2 < \frac12+p$}
  \\
  N_0^\varepsilon &\text{$b_1 < \frac12$, $b_1+b_2 < \frac12$, $b_0+b_1+b_2 = \frac12+p$}
  \\
  1 &\text{$b_1 < \frac12$, $b_1+b_2 < \frac12$, $b_0+b_1+b_2 > \frac12+p$}
  \\
  r^{-\varepsilon}N_0^{p-b_0} &\text{$b_1 < \frac12$, $b_1+b_2 = \frac12$, $b_0 < p$}
  \\
  N_0^\varepsilon &\text{$b_1 < \frac12$, $b_1+b_2 = \frac12$, $b_0 = p$}
  \\
  1 &\text{$b_1 < \frac12$, $b_1+b_2 = \frac12$, $b_0 > p$}
  \\
  r^{\frac12-b_1-b_2} N_0^{\frac12+p-b_0-b_1-b_2} &\text{$b_1 < \frac12$, $b_1+b_2 >\frac12$, $b_0+b_1+b_2 < \frac12+p$}
  \\
  r^{b_0-p+\varepsilon} N_0^\varepsilon &\text{$b_1 < \frac12$, $b_1+b_2 >\frac12$, $b_0+b_1+b_2 = \frac12+p$}
  \\
  r^{b_0-p} &\text{$b_1 < \frac12$, $b_1+b_2 >\frac12$, $\frac12+p-b_1-b_2 < b_0 < p$}
  \\
  r^{-\varepsilon} &\text{$b_1 < \frac12$, $b_1+b_2 >\frac12$, $b_0 = p$}
  \\
  1 &\text{$b_1 < \frac12$, $b_1+b_2 >\frac12$, $b_0 > p$}
  \\
  r^{-b_2+\varepsilon} N_0^{p-b_0-b_2+\varepsilon} &\text{$b_1 = \frac12$, $b_0+b_2 \le p$}
  \\
  r^{b_0-p} &\text{$b_1 = \frac12$, $p-b_2 < b_0 < p$}
  \\
  r^{-\varepsilon} &\text{$b_1 = \frac12$, $b_0 = p$}
  \\
  1 &\text{$b_1 = \frac12$, $b_0 > p$}
  \\
  r^{-b_2} N_0^{p-b_0-b_2} &\text{$b_1 > \frac12$, $b_0+b_2 < p$}
  \\
  r^{b_0-p+\varepsilon} N_0^\varepsilon &\text{$b_1 > \frac12$, $b_0+b_2 = p$}
  \\
  r^{b_0-p} &\text{$b_1 > \frac12$, $p-b_2 < b_0 < p$}
  \\
  r^{-\varepsilon} &\text{$b_1 > \frac12$, $b_0 = p$}
  \\
  1 &\text{$b_1 > \frac12$, $b_0 > p$},
  \end{cases}
$$
for any $\varepsilon > 0$. Plugging this into~\eqref{P:470}, with $r=N_0/N_1$ and $p=\frac12$, we get~\eqref{P:266} for some $A,B \in \R$, so we check that (i) $B \ge A$, (ii) $2B \ge A$ and (iii) we exclude $A=B=0$. Note that (i) implies (ii) if $B \ge 0$. In particular, (i) implies (ii) if $A \ge 0$, since then $B \ge 0$.

\begin{itemize}
\item
If $\rho_{\frac12} \lesssim 1,$ then $A=1-s_0$ and $B=s_1+s_2 \ge 0$, so $B \ge A$ is~\eqref{P:406}. Moreover, $B=0$ implies $A < 0$, by~\eqref{P:423}.
\item
If $\rho_{\frac12} \lesssim r^{-\varepsilon},$ then $A=1-s_0-\varepsilon$ and $B=s_1+s_2-\varepsilon$, so $B \ge A$ is~\eqref{P:406}. Moreover, $2B > A$ by \eqref{P:474}.
\item
If $\rho_{\frac12} \lesssim N_0^\varepsilon,$ then $A=1-s_0+\varepsilon$ and $B=s_1+s_2 \ge 0$. But now $b_0+b_1+b_2 = 1$, so~\eqref{P:419} implies $B > A$ for $\varepsilon > 0$ small enough.
\item
If $\rho_{\frac12} \lesssim N_0^{1-b_0-b_1-b_2},$ then $A=2-s_0-b_0-b_1-b_2$ and $B=s_1+s_2 \ge 0$, so $B \ge A$ is~\eqref{P:402}. Moreover, $B=0$ implies $A < 0$, by~\eqref{P:423}.
\item
If $\rho_{\frac12} \lesssim r^{-\varepsilon} N_0^{\frac12-b_0},$ then $A=\frac32-s_0-b_0-\varepsilon$ and $B=s_1+s_2-\varepsilon$. Now $b_1+b_2=\frac12$, so~\eqref{P:402} implies $B \ge A$, and~\eqref{P:476} implies $2B > A$.
\item
If $\rho_{\frac12} \lesssim r^{\frac12-b_1-b_2} N_0^{1-b_0-b_1-b_2},$ then $A=\frac52-s_0-b_0-2b_1-2b_2$ and $B=s_1+s_2+\frac12-b_1-b_2$, so $B \ge A$ and $2B \ge A$ are the same as~\eqref{P:402} and~\eqref{P:407}. Moreover, $A=0$ implies $B > 0$, in view of~\eqref{P:420}.
\item
If $\rho_{\frac12} \lesssim r^{b_0-\frac12},$ then $A=\frac12-s_0+b_0$ and $B=s_1+s_2+b_0-\frac12$, so $B \ge A$ and $2B \ge A$ are the same as~\eqref{P:406} and~\eqref{P:407}. Moreover, $A=0$ implies $B > 0$, by~\eqref{P:420}.
\item
If $\rho_{\frac12} \lesssim r^{b_0-\frac12+\varepsilon} N_0^\varepsilon,$ then $A=\frac12-s_0+b_0+2\varepsilon$ and $B=s_1+s_2+b_0-\frac12+\varepsilon$, so $2B \ge A$  is~\eqref{P:407}. Now $b_0+b_2=\frac12$ or $b_0+b_1+b_2 = 1$, and in either case~\eqref{P:406} is strict, by~\eqref{P:417} and ~\eqref{P:419}, so $B > A$ for $\varepsilon > 0$ small enough.
\item
If $\rho_{\frac12} \lesssim r^{-b_2} N_0^{\frac12-b_0-b_2},$ then $A=\frac32-s_0-b_0-2b_2$ and $B=s_1+s_2-b_2$, so $B \ge A$ and $2B \ge A$ are the same as~\eqref{P:404} and~\eqref{P:407}, and $A=0$ implies $B > 0$, by~\eqref{P:420}.
\item
If $\rho_{\frac12} \lesssim r^{-b_2+\varepsilon} N_0^{\frac12-b_0-b_2+\varepsilon},$ then $A=\frac32-s_0-b_0-2b_2+2\varepsilon$ and $B=s_1+s_2-b_2+\varepsilon$, so $2B \ge A$ is the same as~\eqref{P:407}. Since $b_1=\frac12$, we infer from~\eqref{P:414} that \eqref{P:404} is strict, hence $B > A$ for $\varepsilon > 0$ small enough.
\end{itemize}

This concludes the proof of Theorem~\ref{P:Thm3}.

\section{The case $b_0 < 0 < b_1,b_2$}\label{P:500:0}

Then the product law reads:

\begin{theorem}\label{P:Thm4} Let $n = 3$. Assume
\begin{align}
  \label{P:500}
  &b_0 < 0 < b_1, b_2
  \\
  \label{P:501}
  & b_0 + b_1 + b_2 \ge \frac12
  \\
  \label{P:502}
  & b_0 + b_1 \ge 0
  \\
  \label{P:503}
  & b_0 + b_2 \ge 0
  \\
  \label{P:504}
  & s_0 + s_1 + s_2  \ge 2 - (b_0 + b_1 + b_2)
  \\
  \label{P:505}
  & s_0 + s_1 + s_2 \ge \frac32 - (b_0 + b_1)
  \\
  \label{P:506}
  & s_0 + s_1 + s_2 \ge \frac32 - (b_0 + b_2)
  \\
  \label{P:507}
  & s_0 + s_1 + s_2 \ge 1 - b_0
  \\
  \label{P:509}
  & s_0 + 2(s_1+ s_2) + b_0 \ge \frac32& &\Hproduct{\text{L} & \text{H} & \text{H} \\ & + & -}
  \\
  \label{P:510}
  & s_1 + s_2 \ge -b_0& &\Hproduct{\text{L} & \text{H} & \text{H} \\ & + & +}
  \\
  \label{P:511}
  & s_0 + s_2 \ge 0& &\Hproduct{\text{H} & \text{L} & \text{H} \\ \phantom{+} & & \phantom{-}}
  \\
  \label{P:512}
  & s_0 + s_1 \ge 0& &\Hproduct{\text{H} & \text{H} & \text{L} \\ \phantom{+} & \phantom{-} &},
\end{align}
as well as the exceptions:
\begin{align}
  \label{P:513}
  &\parbox[t]{10.8cm}{If $b_1=\frac12$, then~\eqref{P:504}=\eqref{P:506} and \eqref{P:505}=\eqref{P:507} must be strict.}
  \\
  \label{P:514}
  &\parbox[t]{10.8cm}{If $b_2=\frac12$, then~\eqref{P:504}=\eqref{P:505} and~\eqref{P:506}=\eqref{P:507} must be strict.}
  \\
  \label{P:517}
  &\parbox[t]{10.8cm}{If $b_1+b_2=1$, then~\eqref{P:504}=\eqref{P:507} must be strict.}
  \\
  \label{P:519}
  &\parbox[t]{10.8cm}{We require~\eqref{P:509} to be strict if $s_0+b_0$ takes one of the values $\frac12$, $\frac32 - 2b_1$, $\frac32 - 2b_2$ or $\frac52 - 2(b_1+b_2)$.}
  \\
  \label{P:520}
  &\parbox[t]{10.8cm}{If~\eqref{P:501} is an equality, then~\eqref{P:502} and~\eqref{P:503} must be strict.}
  \\
  \label{P:521}
  &\parbox[t]{10.8cm}{If one of~\eqref{P:504}--\eqref{P:507} is an equality, then~\eqref{P:510}--\eqref{P:512} must be strict.}
\end{align}
Then $\Hproduct{s_0 & s_1 & s_2 \\ b_0 & b_1 & b_2}$ is a product.
\end{theorem}

It turns out that we can relax the hypotheses somewhat when
\begin{equation}\label{P:530}
  b_0 < 0,
  \qquad
  b_1,b_2 > \frac12.
\end{equation}
Then we first note the following:
\begin{itemize}
\item
\eqref{P:507} implies that~\eqref{P:504}--\eqref{P:506} are strict, and \eqref{P:521} simplifies accordingly.
\item
The exceptional values $\frac32 - 2b_1$, $\frac32 - 2b_2$ and $\frac52 - 2(b_1+b_2)$ of $s_0+b_0$ from~\eqref{P:519} are all strictly less than $\frac12$, so they imply that~\eqref{P:509} is strict, since
$$
  s_0+2(s_1+s_2)+b_0
  \ge
  1+s_1+s_2
  \ge
  2 - (s_0+b_0),
$$
where we applied ~\eqref{P:507} twice. Therefore,~\eqref{P:519} simply says that we must avoid the combination $s_0+b_0=\frac12$ and $s_1+s_2=\frac12$.
\end{itemize}

So \eqref{P:521} and \eqref{P:519} simplify under the assumption \eqref{P:530}. But the following improved result, Theorem~\ref{P:Thm5}, says that we can in fact relax~\eqref{P:521} to:
\begin{equation}\label{P:521:reduced}
  \parbox[t]{10.8cm}{If~\eqref{P:507} is an equality, then~\eqref{P:510} must be strict.}
\end{equation}
Moreover, we can completely ignore~\eqref{P:519}. That is, we can allow the combination
$s_0+b_0=\frac12$ and $s_1+s_2=\frac12$. Note, incidentally, that this implies equality in \eqref{P:507}.

Thus, we claim the following:

\begin{theorem}\label{P:Thm5} If~\eqref{P:530} holds, then the conclusion of Theorem \ref{P:Thm4} remains valid even if we relax its hypotheses as follows: We can dispose of the assumption~\eqref{P:519}, and~\eqref{P:521} can be relaxed to~\eqref{P:521:reduced}.
\end{theorem}

\subsection{Proof of Theorem \ref{P:Thm4}}\label{P:531}

By \eqref{N:92} and \eqref{N:96},
\begin{equation}\label{P:532}
  L_0 \lesssim L_1 + L_2
  +
  \begin{cases}
  \Nmin^{12} &\text{if $\pm_1=\pm_2$}
  \\
  \Nmin^{012} &\text{if $\pm_1\neq\pm_2$},
  \end{cases}
\end{equation}
hence
$P = \Hproduct{s_0 & s_1 & s_2 \\ b_0 & b_1 & b_2}$
is a product if the following are:
\begin{align*}
  P_1 &= \LHproduct{s_0 & s_1 & s_2 \\ 0 & b_0 + b_1 & b_2}
  \\
  P_2 &= \LHproduct{s_0 & s_1 & s_2 \\ 0 & b_1 & b_0 + b_2}
  \\
  P_3 &= \left.\LHproduct{s_0 & s_1+b_0 & s_2 \\ 0 & b_1 & b_2}\right\vert^{(+,+)}_{N_0 \ll N_1 \sim N_2}
  \\
  P_4 &= \left.\LHproduct{s_0+b_0 & s_1 & s_2 \\ 0 & b_1 & b_2}\right\vert_{\text{LHH}}
  \\
  P_5 &= \left.\LHproduct{s_0 & s_1+b_0 & s_2 \\ 0 & b_1 & b_2}\right\vert_{\text{HLH}}
  \\
  P_6 &= \left.\LHproduct{s_0 & s_1 & s_2+b_0 \\ 0 & b_1 & b_2}\right\vert_{\text{HHL}}.
\end{align*}

By symmetry it suffices to consider $P_1$, $P_3$, $P_4$ and $P_5$.

Of course, $P$ is assumed to satisfy the hypotheses of Theorem \ref{P:Thm4}, and we go ahead and check if $P_1$, $P_3$, $P_4$ and $P_5$ satisfy the hypotheses of Theorem~\ref{P:Thm2}. Keeping in mind that $b_0 < 0$, this is readily seen to be the case for $P_1$ (if $b_0 + b_1 > 0$, to be precise; if $b_0+b_1 = 0$, then we use instead Theorem \ref{P:Thm1}). Thus, $P_1$ is a product.

For $P_3$,~\eqref{P:207} and~\eqref{P:210} may fail, but these are not needed in the interaction $N_0 \ll N_1 \sim N_2$ with equal signs. So $P_3$ is a product.

For $P_4$,~\eqref{P:209} and~\eqref{P:210} may fail, but they are not needed since we assume the LHH interaction. So $P_4$ is a product.

For $P_5$,~\eqref{P:207} and~\eqref{P:210} may fail, but they are not needed since we assume the HLH interaction. So $P_5$ is a product.

This concludes the proof of Theorem \ref{P:Thm4}.

\subsection{Proof of Theorem \ref{P:Thm5}}\label{P:540}

Assume \eqref{P:530}. In view of the remarks preceding the theorem, we may assume equality in \eqref{P:507},
$$
  s_0+b_0+s_1+s_2=1,
$$
as otherwise Theorem \ref{P:Thm5} reduces to the already proved Theorem \ref{P:Thm4}.

There are then two things that remain to be proved:
\begin{itemize}
  \item We can allow $s_0+s_2=0$ (hence also $s_0+s_1=0$, by symmetry).
  \item We can allow the combination $s_0+b_0=1/2$ and $s_1+s_2=1/2$.
\end{itemize}

\subsubsection{The case $s_0+s_2=0$} We assume the HLH case, $N_1 \le N_0 \sim N_2$, since this is where \eqref{P:511} is needed. In fact, we can assume
$$
  N_1 \ll N_0 \sim N_2,
$$
since if $N_1 \sim N_0 \sim N_2$, then the LHH case applies, and \eqref{P:511} then plays no role.
 
By \eqref{P:532},
\begin{equation}\label{P:550}
  L_0
  \lesssim
  L_1 + L_2 + N_1,
\end{equation}
so if $\Lmax^{12} \gtrsim N_1$, we reduce (here we rely on \eqref{P:502} and \eqref{P:503}) to checking that $\left.\Hproduct{s_0 & s_1+b_0+b_1 & s_2 \\ 0 & 0 & b_2}\right\rvert_{\text{HLH}}$ and $\left.\Hproduct{s_0 & s_1+b_0+b_2 & s_2 \\ 0 & b_1 & 0}\right\rvert_{\text{HLH}}$ are products. But this follows from Theorem \ref{P:Thm1}, since $s_0 + s_1 + s_2 + b_0 + b_j = 1 + b_j > \frac32$ for $j=1,2$.

Thus, it remains to consider the regime
\begin{equation}\label{P:552}
  \Lmax^{12} \ll N_1 \ll N_0 \sim N_2.
\end{equation}
Then by \eqref{P:550} we could reduce to proving that $\left.\Hproduct{s_0 & s_1+b_0 & s_2 \\ 0 & b_1 & b_2}\right\rvert_{\text{HLH}}$ is a product, but this approach fails, since $s_0+s_2=0$.

To see what goes wrong, let us recall our usual method. We want to prove the estimate
\begin{equation}\label{P:554}
  \norm{uv}_{H^{-s_0,-b_0}} \le C \norm{u}_{H^{s_1,b_1}}\norm{v}_{H^{s_2,b_2}}.
\end{equation}
If we apply \eqref{P:550} and then follow our usual approach of writing the $L^2$ product estimate as a trilinear integral estimate by duality, and then apply the dyadic estimates and try to sum the pieces, we come up short. In fact, we will be left with the sum (since we are in the HLH case, and since $s_0+s_2=0$ and $s_1+b_0=1$)
$$
  \sum_{\Nbold} \chi_{N_1 \le N_0 \sim N_2}
  \bignorm{F_0^{N_0}}
  \bignorm{F_1^{N_1}}
  \bignorm{F_2^{N_2}}.
$$
But of course then we have no way of summing $N_1$.

To avoid this problem, we delay the application of \eqref{P:550}, and begin instead by writing \eqref{P:554} as a doubled estimate:
$$
  \abs{I}
  \lesssim
  \norm{u}_{H^{s_1,b_1}}^2\norm{v}_{H^{s_2,b_2}}^2
$$
where
\begin{equation}\label{P:556}
  I
  =
  \iint \angles{D}^{-s_0}\angles{D_-}^{-b_0}(uv)
  \cdot
  \overline{\angles{D}^{-s_0}\angles{D_-}^{-b_0}(uv)}
  \d t \d x.
\end{equation}
The crucial point now is that by Plancherel we can move the multiplier $\angles{D_-}^{-b_0}$ from the second product onto the first:
\begin{equation}\label{P:558}
  I
  =
  \iint \angles{D}^{-s_0}\angles{D_-}^{-2b_0}(uv)
  \cdot
  \overline{\angles{D}^{-s_0}(uv)}
  \d t \d x,
\end{equation}
and vice versa.

Now we make the dyadic decomposition for both products, restricted by~\eqref{P:552}. Let us denote the dyadic sizes by $\Nbold=(N_0,N_1,N_2)$, $\Lbold=(L_1,L_2)$ for the leftmost product in $I$ and $\Nbold'=(N_0',N_1',N_2')$ and $\Lbold'=(L_1',L_2')$ for the rightmost product. By Plancherel, we must have $N_0=N_0'$, and by symmetry we may assume $N_1 \le N_1'$.

Then we can see the advantage of writing $I$ as in \eqref{P:558}: Since we now have
$$
  \Lmax^{12} \ll N_1 \le N_1',
$$
the symbol of $\angles{D_-}^{-2b_0}$ in \eqref{P:558} will have a size $O(N_1^{-2b_0})$, whereas if we had kept $I$ in the form \eqref{P:556}, then we would have had two instances of the multiplier $\angles{D_-}^{-b_0}$, with symbol sizes $O(N_1^{-b_0})$ and $O((N_1')^{-b_0})$, respectively. Thus, we have essentially gained a factor $(N_1/N_1')^{-b_0}$, and this makes it possible to sum without running into any divergences, as we now show.

In fact, after the dyadic decomposition, we are faced with a sum
$$
  S = \sum_{\Nbold,\Lbold,\Lbold'}
  \chi_{N_1 \le N_1' \ll N_0 \sim N_2 \sim N_2'}
  N_1^{-2b_0} N_0^{-2s_0} \Abs{I_{\Nbold,\Lbold,\Lbold'}},
$$
where $\Nbold=(N_0,N_1,N_1',N_2,N_2')$ and
$$
  I_{\Nbold,\Lbold,\Lbold'}
  =
  \iint u_1^{N_1,L_1}u_2^{N_2,L_2}
  \cdot
  \overline{u_1^{N_1',L_1'}u_2^{N_2',L_2'}}
  \d t \d x.
$$
Here $u^{N,L}$ is defined as in \S\ref{N:1}.

By the wave type estimate \eqref{N:50}--\eqref{N:54} (rewritten as an $L^2$ bilinear estimate),
\begin{align*}
  \abs{I_{\Nbold,\Lbold,\Lbold'}}
  &\lesssim N_1 N_1' (L_1L_2 L_1'L_2')^{1/2}
  \bignorm{u_1^{N_1,L_1}}
  \bignorm{u_2^{N_2,L_2}}
  \bignorm{u_1^{N_1',L_1'}}
  \bignorm{u_2^{N_2',L_2'}}
  \\
  &\lesssim (N_1 N_1')^{1-s_1} (N_0)^{-2s_2} (L_1L_1')^{1/2-b_1} (L_2L_2')^{1/2-b_2}
  \\
  &\;\;
  \times\bignorm{u_1^{N_1,L_1}}_{H^{s_1,b_1}}
  \bignorm{u_2^{N_2,L_2}}_{H^{s_2,b_2}}
  \bignorm{u_1^{N_1',L_1'}}_{H^{s_1,b_1}}
  \bignorm{u_2^{N_2',L_2'}}_{H^{s_2,b_2}},
\end{align*}
hence, keeping in mind that $s_0+s_2=0$ and $s_1+b_0=1$,
\begin{multline*}
  S \lesssim
  \sum_{\Nbold,\Lbold,\Lbold'}
  \chi_{N_1 \le N_1' \ll N_0 \sim N_2 \sim N_2'}
  \left(\frac{N_1}{N_1'}\right)^{-b_0} (L_1L_1')^{1/2-b_1} (L_2L_2')^{1/2-b_2}
  \\
  \times\bignorm{u_1^{N_1,L_1}}_{H^{s_1,b_1}}
  \bignorm{u_2^{N_2,L_2}}_{H^{s_2,b_2}}
  \bignorm{u_1^{N_1',L_1'}}_{H^{s_1,b_1}}
  \bignorm{u_2^{N_2',L_2'}}_{H^{s_2,b_2}}.
\end{multline*}
Since $b_1,b_2 > 1/2$, we can trivially sum the $L$'s, and $N_0 \sim N_2 \sim N_2'$ can be summed by Cauchy-Schwarz. Thus,
$$
  S \lesssim
  \bignorm{u_2}_{H^{s_2,b_2}}^2
  \sum_{N_1, N_1'}
  \chi_{N_1 \le N_1'}
  \left(\frac{N_1}{N_1'}\right)^{-b_0}
  \bignorm{u_1^{N_1}}_{H^{s_1,b_1}}
  \bignorm{u_1^{N_1'}}_{H^{s_1,b_1}},
$$
and writing $N_1'=MN_1$, where $M \ge 1$ is dyadic,
$$
  S \lesssim
  \bignorm{u_2}_{H^{s_2,b_2}}^2
  \sum_{N_1, M}
  \left(\frac{1}{M}\right)^{-b_0}
  \bignorm{u_1^{N_1}}_{H^{s_1,b_1}}
  \bignorm{u_1^{M N_1}}_{H^{s_1,b_1}}.
$$
Now apply Cauchy-Schwarz in $N_1$, and observe that $\sum_{M} \left(\frac{1}{M}\right)^{-b_0}$ converges.

This concludes the proof that we can allow $s_0+s_2=0$.

\subsubsection{The combination $s_0+b_0=\frac12$ and $s_1+s_2=\frac12$} The issue is that we have equality in \eqref{P:509}, so we restrict to the LHH interaction with opposite signs. In particular, we may as well assume $s_1=s_2=\frac14$.

We apply the hyperbolic Leibniz rule, \eqref{N:92}. Corresponding to the first two terms in its right hand side, we need to check that $\left.\Hproduct{s_0 & s_1 & s_2 \\ 0 & b_0+b_1 & b_2}\right\rvert_{\text{LHH}}$ and $\left.\Hproduct{s_0 & s_1 & s_2 \\ 0 & b_1 & b_0+b_2}\right\rvert_{\text{LHH}}$ are products, but this follows from Theorem \ref{P:Thm2} (or Theorem \ref{P:Thm1} if $b_0+b_1=0$ or $b_0+b_2=0$).

We are left with the third term in the right hand side of \eqref{N:92}, corresponding to which we define bilinear operators as in \eqref{N:97}. Thus, we need to show
$$
  \norm{\mathfrak B_{(+,-)}^{-b_0}(u,v)}_{H^{-s_0,0}}
  \lesssim \norm{u}_{H^{1/4,b_1}} \norm{v}_{H^{1/4,b_2}}
$$
for $u,v$ such that $\widetilde u(\tau,\xi)$ and $\widetilde v(\tau,\xi)$ are supported in $\tau \ge 0$ and $\tau \le 0$, respectively.
 
Since $s_0$ is positive, it suffices to show the homogeneous variant
$$
  \norm{D^{-s_0}\mathfrak B_{(+,-)}^{-b_0}(u,v)}_{L^2(\R^{1+3})}
  \lesssim \bignorm{D^{1/4} u}_{H^{0,b_1}} \bignorm{D^{1/4} v}_{H^{0,b_2}},
$$
where $D^\alpha$ is the multiplier corresponding to the symbol $\abs{\xi}^\alpha$, for $\alpha \in \R$. By the transfer principle, this follows from the corresponding estimate for two solutions of the homogeneous wave equation:
$$
  \norm{D^{-s_0}\mathfrak B_{(+,-)}^{-b_0}(u^+,v^-)}_{L^2(\R^{1+3})}
  \lesssim \bignorm{D^{1/4} f}_{L^2(\R^3)} \bignorm{D^{1/4} g}_{L^2(\R^3)},
$$
where $u^+(t) = e^{it\abs{\nabla}}f$ and $v^-(t) = e^{-it\abs{\nabla}}g$. This last estimate is proved in \cite{Foschi:2000}; it relies on the assumptions $b_0 < 0$ and $s_0+b_0=\frac12$.

\section{Reformulation of the rules on the boundary}\label{Z}

Here we prove Theorem \ref{B:Thm3}.

We first show that the rules~\eqref{B:53},~\eqref{B:58} and~\eqref{B:65} imply~\eqref{B:80}--\eqref{B:96}.

To prove~\eqref{B:80}, we assume $b_0=\frac12$, and show that \eqref{B:1}, \eqref{B:5}, \eqref{B:6} and \eqref{B:7} must then be strict. To this end, we assume each in turn to be an equality, and deduce a contradiction. First, if \eqref{B:1} is an equality, then~\eqref{B:53} implies that~\eqref{B:52} is strict, that is, $\frac12 > \max(\frac12,b_1,b_2)$, which is impossible. Second, if \eqref{B:5} is an equality, then so are~\eqref{B:54} and ~\eqref{B:60}, hence the rules~\eqref{B:58} and~\eqref{B:65} imply that~\eqref{B:55} and~\eqref{B:61} are strict, but this contradicts equality in~\eqref{B:5}. Third, if~\eqref{B:6} is an equality, then so are~\eqref{B:55},~\eqref{B:56},~\eqref{B:61} and~\eqref{B:62}, violating the rules~\eqref{B:58} and~\eqref{B:65}. A similar argument shows that~\eqref{B:7} is strict.

Thus we have proved~\eqref{B:80}, and~\eqref{B:81}--\eqref{B:89} follow by similar arguments which we leave to the interested reader.

The exceptions~\eqref{B:90}--\eqref{B:95} arise in a similar way when we compare~\eqref{B:63} with~\eqref{B:60}--\eqref{B:62} and~\eqref{B:64}. The thing to note here, however, is that such a comparison gives not only the exceptional values of $s_0-b_0$ listed in ~\eqref{B:90}--\eqref{B:95}, but also the values $\frac{n}2 - 2(b_1+b_2)$, $\frac{n-2}2 - 2b_1$ and $\frac{n-2}2 - 2b_2$, which imply that~\eqref{B:13}, or equivalently \eqref{B:63}, coincides with~\eqref{B:8},~\eqref{B:10} and~\eqref{B:11}, respectively. We claim, however, that these values automatically imply strict inequality in~\eqref{B:13}, hence there is no need to list them.

Consider first $s_0-b_0=\frac{n}2 - 2(b_1+b_2)$, so that ~\eqref{B:13} coincides with~\eqref{B:8}. To get a contradiction, assume that they are equalities. Then \eqref{B:5} must be strict, since otherwise we would be in the case \eqref{B:90}, but then we know \eqref{B:13} cannot be an equality. So \eqref{B:8} is an equality and \eqref{B:5} is strict, hence $b_0 > \frac12$. Similarly, comparing \eqref{B:13}, or equivalently \eqref{B:63}, with \eqref{B:12}, we get $\frac{s_0-b_0}2 > \frac14$. Thus, $s_0+b_0 > 1 + \frac12 = \frac32$. But this implies that
$$
  (s_0+b_0)+2(s_1+s_2) \ge s_0+b_0 > \frac32 \ge \frac{n}2,
$$
where the first inequality holds by \eqref{B:19}, and the last by the fact that $n \le 3$. This shows that \eqref{B:13} is strict, so we have a contradiction.

A similar argument, which we omit, shows that \eqref{B:13} must be strict if $s_0-b_0$ takes one of the values $\frac{n-2}2 - 2b_1$ or $\frac{n-2}2 - 2b_2$.

Finally, consider \eqref{B:96}. If one of~\eqref{B:5}--\eqref{B:12} is an equality, then so is one of~\eqref{B:54}--\eqref{B:56} and one of \eqref{B:60}--\eqref{B:62}, so the rules~\eqref{B:58} and~\eqref{B:65} guarantee that \eqref{B:57} and \eqref{B:64} are strict, hence so are \eqref{B:16} and \eqref{B:19}.

Conversely, we must prove that \eqref{B:80}--\eqref{B:96} imply the rules~\eqref{B:53},~\eqref{B:58} and~\eqref{B:65}.

First assume that equality holds in both \eqref{B:51} and \eqref{B:52}, to get a contradiction. By permutation, it suffices to consider the case $\max(b_0,b_1,b_2) = b_0$. Then it follows that $b_0 = \frac12$ and $b_1+b_2=0$, but this contradicts \eqref{B:80}.

We proceed similarly to prove \eqref{B:58} and \eqref{B:65}. Let us just show one representative example. Say we have equality in both \eqref{B:55} and \eqref{B:56}. Then we will have equality in two of \eqref{B:5}--\eqref{B:12}. For example, if the maximum in the right hand sides of \eqref{B:55} and \eqref{B:56} are $-b_0-b_1$ and $-\frac{n-3}{4}$, respectively, then we have equality in \eqref{B:6} and \eqref{B:12}. But this means that $b_0+b_1=\frac{n-1}{4}$, so \eqref{B:86} implies that \eqref{B:6} and \eqref{B:12} are strict, and we have a contradiction.

This concludes the proof of Theorem \ref{B:Thm3}.

\bibliographystyle{amsplain}
\bibliography{atlasbibliography}

\end{document}